\documentclass[12pt]{amsart}
\usepackage{amsmath}
\usepackage{amsfonts}
\usepackage{amssymb}
\usepackage{mathrsfs}
\usepackage{amsthm}
\usepackage{comment}
\usepackage{chngcntr}
\counterwithout{equation}{section}
\usepackage{graphicx}
\usepackage[svgnames,dvipsnames]{xcolor}
\usepackage{float}
\usepackage{diagbox}  
\usepackage[textwidth=22mm]{todonotes}
\usepackage[normalem]{ulem}
\usepackage{hyperref}
\usepackage{mathtools}

\DeclarePairedDelimiter\floor{\lfloor}{\rfloor}

\hypersetup{
    linktocpage=false,
    colorlinks=true,
    linkcolor={BrickRed},
    citecolor={ForestGreen!90!black},
    urlcolor={Salmon!50!black}
}

\numberwithin{equation}{section}

\newtheorem{proposition}{Proposition}
\newtheorem{theorem}{Theorem}
\newtheorem{lemma}{Lemma}

\theoremstyle{definition}
\newtheorem{remark}{Remark}

\DeclareMathOperator{\tr}{tr}

\newcommand{\ad}{\mathrm{ad}}

\newcommand{\adstar}{\ad^{\star}}

\newcommand\conn[2]{\Gamma(#1, #2)}

\newcommand{\Z}{\mathbb{Z}}

\newlength{\storeparskip}
\title{The Large-$N$ Limit of Ricci Curvature in the Zeitlin Model on the Torus}

\author{Leandro Lichtenfelz}
\author{Isabelle Raad}
\author{Justin Valletta} 

\address{L. Lichtenfelz: Wake Forest University, Winston-Salem, NC., U.S.A.}
\email{lichtel@wfu.edu}
\address{I. Raad: Wake Forest University, Winston-Salem, NC., U.S.A.}
\email{raadiz25@wfu.edu}
\address{J. Valletta: Wake Forest University, Winston-Salem, NC., U.S.A.}
\email{valletj@wfu.edu}
 
\begin{document}  

\begin{abstract}

The Euler equations of hydrodynamics describe the motion of an ideal incompressible fluid. Arnold observed that their solutions are geodesics on the group of volume-preserving diffeomorphisms of the fluid domain, equipped with the $L^2$ metric, so that the curvature of this group governs how nearby fluid motions separate. Sectional curvature takes both signs, and Ricci curvature, its average over all directions, is the natural measure of which sign prevails. In infinite dimensions, however, the sum defining it diverges. When the fluid domain is the flat torus, we define Ricci curvature through Zeitlin's finite-dimensional models, which approximate the Euler equations by geodesic equations on the matrix groups $SU(N)$. We compute the Ricci curvature of these models exactly and prove that, normalized by dimension, it converges as $N \to \infty$ to an explicit integral. The limit is strictly negative in every Fourier mode and decreases logarithmically with the wavenumber, the behavior conjectured in \cite{lichtenfelz2026ricci} for the sphere. We also show that the diameter of $SU(N)$ in the Zeitlin metric grows at most linearly in $N$, which bounds the time needed to reach any configuration. \\[2ex]

    \\[2ex]
\end{abstract}

\maketitle 

\tableofcontents

\newpage

\section{Introduction}

The incompressible Euler equations describe the motion of an ideal fluid on a domain $M$, taken here to be a compact, oriented Riemannian manifold. If $u$ denotes the time-dependent vector field on $M$ representing the velocity of the fluid, the Euler equations are
\begin{gather}\label{eq_euler}
\begin{split}
\partial_tu + \nabla_uu &= -\nabla p \\
\mathrm{div}\, u &= 0 \\
u(0) &= u_0,
\end{split}
\end{gather}
where $\nabla_uu$ is the covariant derivative of $u$ along itself, $p$ is the pressure function, and $u_0$ is the initial velocity field.

Arnold~\cite{arnold2013differential} reinterpreted the Euler equations~\eqref{eq_euler} as the
geodesic flow on the infinite-dimensional Lie group $\mathrm{Diff}_{\mu}(M)$ equipped with the right-invariant $L^2$ metric. Curvature influences the behavior of nearby geodesics through the Jacobi equation, providing a geometric perspective on the stability of fluid flows. This connection led him to compute sectional curvatures on $\mathrm{Diff}_{\mu}(T^2)$, followed by many others who worked on a broader class of manifolds (see, e.g., \cite{lukatskii1984curvature,rouchon1992thejacobi,misiolek1993stability,preston2005nonpositive,smolentsev2006diffeomorphism}).

From these computations, a recurring theme emerges: while sectional curvature on $\mathrm{Diff}_{\mu}(M)$ takes on both signs, it tends to be negative. Ricci curvature offers a natural way to quantify this tendency by summing sectional curvatures over a complete set of orthonormal directions. However, this sum generally diverges in the infinite-dimensional setting, requiring an appropriate renormalization.

An early renormalization was obtained by Lukatskii \cite{lukatskii1984curvature}, who averaged Arnold's sectional curvatures over growing sets of Fourier modes on the flat torus. More recently, Khesin, Lang\o en and Markina \cite{khesin2025curvatures} extended these computations to the Klein bottle and the projective plane. In these constructions, sectional curvatures of the infinite-dimensional diffeomorphism group are averaged over finite-dimensional Laplacian eigenspaces, which are generally not Lie subalgebras.

A different approach, developed for the sphere $S^2$ by the first author together with Modin and Preston \cite{lichtenfelz2026ricci}, builds on Zeitlin's matrix model of hydrodynamics the sphere \cite{Ze2004}. This model equips $SU(N)$ with a right-invariant Riemannian metric whose geodesic equation is a matrix approximation of the Euler equations \eqref{eq_euler}. It therefore has a finite-dimensional Ricci tensor for each $N$, and normalizing by the dimension $N^2-1$ of $SU(N)$ defines what we call the \emph{averaged Ricci curvature}. 

On the sphere, the Ricci tensor of the Zeitlin metric was shown \cite{lichtenfelz2026ricci} to be a multiple $r_{\ell}(N)$ of the metric on the $\ell$th eigenspace of the quantized Laplacian. Furthermore, based on numerical evidence and conjectured identities involving Wigner symbols, the authors proposed that the averaged Ricci curvature $\widetilde{r}_{\ell}(N) := \tfrac{r_{\ell}(N)}{N^2-1}$ satisfies:
\begin{equation*}
\widetilde{r}_{\ell}(N) \longrightarrow -\frac{H_{\ell}-1}{2}, \qquad \text{as $N \rightarrow \infty$},
\end{equation*}
where $H_{\ell}$ denotes the $\ell$th harmonic number. The analogous problem for the flat torus, the setting of Zeitlin's original construction \cite{Ze1991}, was left open.

In this paper, we carry out the corresponding program on the flat torus and prove convergence of averaged Ricci curvature to an explicit, strictly negative limit. We also derive conjugate time estimates using techniques from \cite{lebrigant2025geodesics}. Then, we use these estimates to prove that the diameter of $SU(N)$ under the Zeitlin metric grows at most linearly with $N$. To state our results precisely, we briefly describe the torus model, including a particular Fourier basis that is used throughout the paper.

In the torus model, we work with the group $\mathrm{Diff}_{\mu, \mathrm{ex}}(T^2)$ of Hamiltonian diffeomorphisms of the $T^2$, whose Lie algebra consists of exact divergence-free vector fields, i.e., skew-gradients of periodic stream functions. This excludes the two-dimensional subspace of harmonic vector fields.

On $T^2 = \mathbb{R}^2/2\pi\mathbb{Z}^2$ and for $n\in\Z^2\setminus\{0\}$, the stream functions $\cos(n\cdot x)$ and $\sin(n\cdot x)$ are mean-zero eigenfunctions of the Laplacian with eigenvalue $-\|n\|^2$. The Poisson bracket of two such functions, with frequencies $m$ and $n$, is a linear combination of modes with frequencies $m+n$ and $m-n$ and coefficients proportional to $n \star m = n_1m_2 - n_2m_1$. Working over $\mathbb{C}$, the exponentials $e^{in\cdot x}$ are replaced by $N\times N$ complex matrices $J_n$ from which a basis of $\mathfrak{su}(N)$ is obtained by setting
\begin{equation*}
S_n := \frac{1}{2}\Big(J_n - J_{-n}\Big), \qquad C_n := \frac{1}{2i}\Big(J_n + J_{-n}\Big).
\end{equation*}
For odd $N \geq 3$, define the index sets:
\begin{align*}
\Lambda_N &:= \big\{n\in\mathbb{Z}^2\backslash{}\{0\}
:\ |n_1|,|n_2|\leq\tfrac{N-1}{2}\big\}, \\
\Lambda_N^+ &:= \big\{n\in \Lambda_N \,:\, \text{$n_1 > 0$ or $n_1 = 0$ and $n_2 > 0$}
\big\}.
\end{align*}
Then, the set $\big\{C_n, S_n  : n \in \Lambda^+_N \big\}$ is a basis for $\mathfrak{su}(N)$, and the Zeitlin metric $\langle\cdot,\cdot\rangle$ makes this basis orthogonal, with
\begin{equation*}
\|C_n\|^2 = \|S_n\|^2 = \frac{1}{2} \|n\|^2,
\end{equation*}
extended by right-invariance to $SU(N)$. This is the analogue of the $L^2$ metric on divergence-free vector fields (or the $H^1$ metric on stream functions). Its Ricci tensor has the following form.

\begin{theorem}\label{eq_intro_thm_1}
Let $n\in\Lambda_N^+$. The Ricci tensor $\mathrm{Ric}^{(N)}$ of the Zeitlin metric on $SU(N)$ satisfies $\mathrm{Ric}^{(N)}(u,v) = r_n(N)\langle u,v\rangle$ for all $u\in\operatorname{span}\{C_n,S_n\}$ and $v\in\mathfrak{su}(N)$, where
\begin{equation}\label{eq_r_n_intro}
r_n(N) = \sum_{\substack{m\in\Lambda_N\\ m\neq \pm n}} \big(\alpha^{(N)}_{n,m}\big)^2\,\frac{\|n\|^4-\big(\|m+n\|^2-\|m\|^2\big)^2}{\|n\|^2\,\|m\|^2\,\|m+n\|^2},
\end{equation}
and the coefficients
\begin{equation}\label{eq_alpha_intro}
 \alpha^{(N)}_{n, m} = \frac N{2\pi}\sin\!\left(\frac{2\pi(n\star m)}{N}\right)
\end{equation} 
are the scaled Lie algebra structure constants relative to the basis $\{C_n, S_n\}$.
\end{theorem}

\begin{remark} 
Formula~\eqref{eq_r_n_intro} has the same structure as its counterpart on the sphere~\cite[(3.12)]{lichtenfelz2026ricci}, with the Laplace eigenvalues $\ell(\ell+1)$ on $S^2$ replaced by $\|n\|^2$ and the weights built from squared Wigner $6j$ symbols replaced by the squared structure constants~\eqref{eq_alpha_intro}. It is also important to note that in \eqref{eq_r_n_intro}, the symbol $\|\cdot\|$ denotes the standard Euclidean norm, so $\|n\|^2 = n_1^2 + n_2^2$, but with the condition that $n$ is reduced modulo $N$ to its representative in $\Lambda_N$ before the norm is evaluated. This affects terms like $\|m+n\|$. We use this convention throughout the paper.
\end{remark}

For fixed $n$, the summands in \eqref{eq_r_n_intro} are uniformly bounded in $m$ and $N$, while the number of terms is on the order of $N^2$. As in the sphere case, we therefore normalize by the dimension of $SU(N)$ and define the averaged Ricci curvature by
\begin{equation*}
\widetilde{r}_n(N) := \frac{r_n(N)}{N^2-1}.
\end{equation*}
The main result of this paper establishes that $\widetilde{r}_n(N)$ converges to a finite, negative limit as $N\to\infty$, given by an explicit integral formula, which serves as a notion of Ricci curvature on $\mathrm{Diff}_{\mu, \mathrm{ex}}(T^2)$. 

\begin{theorem}\label{eq_main_thm_intro}
For every $n\in\Z^2\setminus\{0\}$, as $N\to\infty$ through odd integers, the averaged Ricci curvature converges
\begin{equation}\label{eq_limit_intro}
\widetilde r_n(N) \longrightarrow -\frac{1}{\pi^2\|n\|^2}\int_{Q} \frac{\sin^2\big(2\pi\,(n\star x)\big)\,(n\cdot x)^2}{\|x\|^4}\,dx, \qquad Q = \big[-\tfrac12,\tfrac12\big]^2.
\end{equation}
In particular, the limit is finite and strictly negative.
\end{theorem}

Since the integral in \eqref{eq_limit_intro} is strictly negative, for each fixed $n$ the averaged Ricci curvature $\widetilde{r}_n(N)$ must be eventually negative. The plot below illustrates this: for a fixed matrix size $N$, curvature is positive in the high-frequency modes, but becomes negative once $N$ is large enough. Numerically, the transition seems to occurs along the line $\|n\| = 0.29N$.

\begin{figure}[H]\label{fig_1}
    \centering
    \includegraphics[scale=0.75]{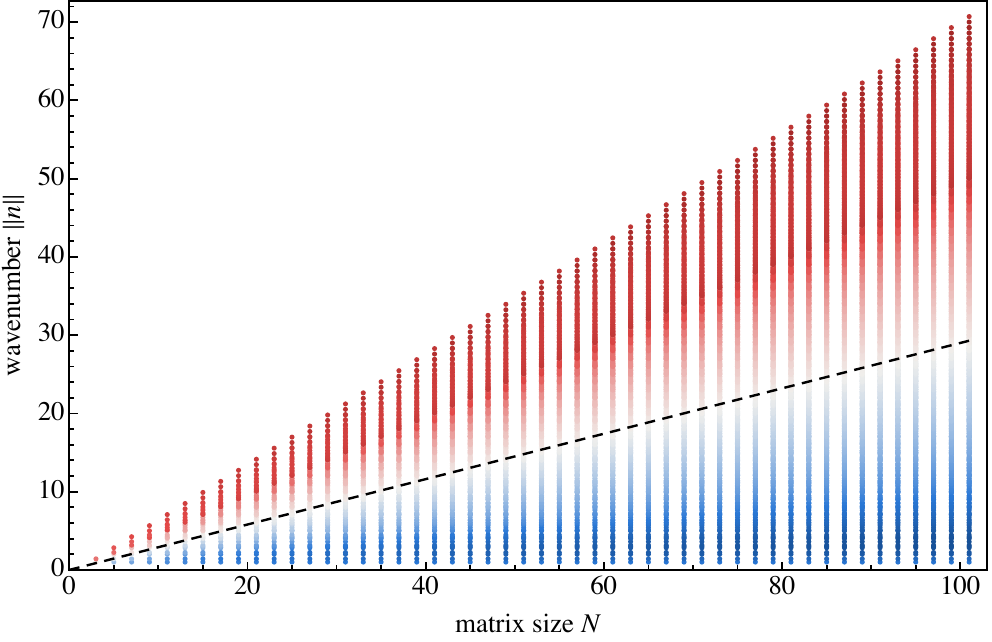}
    \caption{Sign of the averaged Ricci curvature $\widetilde{r}_n(N)$ of the Zeitlin metric on $SU(N)$: positive in red, negative in blue. The dashed line is $\|n\| = 0.29N$.}
\end{figure}

The limiting curvature also reveals a striking similarity between the sphere and the torus models. Despite the different forms of the finite-dimensional Ricci formulas, involving Wigner $6j$ symbols on the sphere and trigonometric functions on the torus, both limits exhibit a common high-frequency pattern. Writing $\rho_{T^2}(n)$ for the limit in
\eqref{eq_limit_intro} and $\rho_{S^2}(\ell)$ for the
conjectured limit on the sphere, their asymptotic expansions after rescaling take the form:
\begin{align}\label{eq_asymp_intro}
2\rho_{S^2}(\ell)
&=-\log\ell-\gamma+1+o(1),\\
2\pi\rho_{T^2}(n)
&=-\log\|n\|-\gamma
+\frac{1}{2}+\frac{2G}{\pi}-\log(2\pi)+o(1),\label{torus_asymp}
\end{align}
as $\ell\to\infty$ and $\|n\|\to\infty$ respectively, where $\gamma$ is the Euler-Mascheroni constant and $G$ is Catalan's constant. We prove \eqref{torus_asymp} in Proposition \ref{prop_asymp}, section $4$. In particular, both expansions contain the same leading logarithmic term.

\subsection*{Acknowledgements}
The authors recognize the concurrent and independent work of
Ishida and Suri~\cite{ishida2026ricci}, of which we originally became aware via private communication. They pursued a similar program by using the Zeitlin model to define a notion of Ricci curvature on the volume-preserving diffeomorphism group of the torus. Their emphasis, however, is on Ricci curvature induced by Sobolev metrics and the incorporation of the Coriolis effect, whereas ours is on the convergence of Ricci curvature and conjugate points. The authors also thank Klas Modin and Stephen C. Preston for helpful discussions. 

\section{Background}

\subsection{Incompressible Euler Equations in 2D}

The two-dimensional Euler equations have a rich geometric structure which is the basis of Zeitlin's construction. We describe it briefly here, and refer the reader to \cite{arnold2021topological}, \cite{MoVi2024}, and \cite{frankel1998geometry} for additional information.

On the flat torus $T^2$ with coordinates $x = (x_1, x_2)$ and volume form $\mu = dx_1 \wedge dx_2$, the Levi-Civita connection term in \eqref{eq_euler} is $\nabla_uu = (u \cdot \nabla)u$, and the scalar vorticity of $u$ is given by
\begin{equation*} 
\omega := \mathrm{curl}\, u = \partial_1 u_2 - \partial_2 u_1.
\end{equation*} 
Taking the curl of the first equation in \eqref{eq_euler} eliminates the pressure, and since 
\begin{equation*}
\mathrm{curl}\big((u\cdot \nabla)u\big) = u\cdot \nabla{\omega} + \omega\, \mathrm{div}\, u
\end{equation*}
in two dimensions, the incompressibility condition gives
\begin{equation}\label{eq_eulerscalar}
\partial_t\omega + u \cdot \nabla \omega = 0.
\end{equation}
Thus, vorticity is transported by the flow. Next, we express the velocity in terms of stream functions. By the Hodge decomposition, a divergence-free vector field $u$ can be written uniquely as $u = \nabla^{\perp} \psi + h$, where $\psi$ is a mean-zero periodic function, $\nabla^\perp\psi := (\partial_2\psi, -\partial_1\psi)$ is the skew-gradient, and $h$ is harmonic, hence constant. We restrict to the subspace $h = 0$, which is invariant under the Euler flow. Its elements are the exact divergence-free vector fields, which form the Lie algebra $T_e\mathrm{Diff}_{\mu, \mathrm{ex}}(T^2)$ of the group $\mathrm{Diff}_{\mu, \mathrm{ex}}(T^2)$ of Hamiltonian diffeomorphisms of $T^2$.

For $u = \nabla^{\perp}\psi$ we have $\omega = -\Delta \psi$, and since $\Delta$ is invertible on mean-zero functions, the vorticity and the stream function determine each other. Finally, let
\begin{equation*} 
\{f, g\} := \partial_1f \partial_2g - \partial_2f \partial_1g
\end{equation*}
denote the Poisson bracket of functions on $T^2$. Then we have
\begin{equation*}
u\cdot \nabla \omega = \partial_2\psi \partial_1 \omega - \partial_1\psi\partial_2\omega = \{ \omega, \psi \},
\end{equation*}
so that \eqref{eq_eulerscalar} becomes
\begin{equation}\label{eq:liepoisson}
\partial_t\omega + \{ \omega, \psi \} = 0, \qquad \omega = -\Delta \psi.
\end{equation}
This is the Euler--Arnold equation of the $L^2$ metric on $\mathrm{Diff}_{\mu, \mathrm{ex}}(T^2)$, written in vorticity form. Note that the Lie algebra $T_e\mathrm{Diff}_{\mu, \mathrm{ex}}(T^2)$ is identified here with mean-zero stream functions via $\psi \mapsto \nabla^{\perp}\psi$, its Lie bracket corresponds to the Poisson bracket, and the metric enters through the relation $\omega = -\Delta \psi$.

\subsection{The Zeitlin model on the torus}
The Zeitlin model for the Euler equations on the torus, originally introduced in~\cite{Ze1991}, replaces the Poisson bracket and the Laplacian by finite-dimensional analogues: first, the Fourier modes are replaced by special matrices whose scaled commutators converge to the Poisson brackets of the corresponding modes, so that the Lie algebra of exact divergence-free vector fields is approximated by $\mathfrak{su}(N)$. This step, in which Poisson brackets are replaced by commutators, can be traced back to Dirac's formulation \cite{dirac1925fundamental} of quantization. 

Second, the Laplacian is replaced by an operator on $\mathfrak{su}(N)$ with eigenvalues $-\|n\|^2$ on the corresponding matrices, which defines the Zeitlin metric on $SU(N)$. The resulting Euler--Arnold equation is a finite-dimensional analogue of~\eqref{eq:liepoisson}, but not its Fourier truncation: the bracket is deformed, and wave vectors are added modulo $N$.

We expand the vorticity and the stream function in Fourier series:
\begin{equation*}
\omega(x,t)=\sum_{m}\omega_m(t)\,e^{im\cdot x}, \qquad \psi(x,t)=\sum_{m}\psi_m(t)\,e^{im\cdot x},    
\end{equation*}
where $m$ ranges over $\mathbb{Z}^2\setminus\{0\}$. The Poisson bracket and the Laplacian act on stream functions by
\begin{equation}\label{eq_pb_fourier}
\{e^{in\cdot x}, e^{im\cdot x}\} = -(n\star m)\,e^{i(n+m)\cdot x}, \qquad -\Delta e^{in\cdot x} = \|n\|^2 e^{in\cdot x}.
\end{equation}
The second relation turns $\omega=-\Delta\psi$ into $\omega_m=\|m\|^2\psi_m$, and with the first,~\eqref{eq:liepoisson} becomes
\begin{equation}\label{eq:fourier}
\partial_t\omega_p=\sum_{k\neq 0,\,p}\frac{p\star k}{\|k\|^2}\,\omega_k\,\omega_{p-k}.
\end{equation}
Zeitlin's model replaces $e^{in\cdot x}$ by the matrices $J_n$ of~\eqref{eq_J_definition}. By~\eqref{eq_lie_bracket},
\begin{equation*}
\frac{N}{4\pi i}\,[J_n, J_m] = -\alpha^{(N)}_{n,m}\,J_{n+m}, \qquad \alpha^{(N)}_{n,m} = \frac{N}{2\pi}\sin\Big(\frac{2\pi(n\star m)}{N}\Big),
\end{equation*}
which is the first relation in~\eqref{eq_pb_fourier}, with $n\star m$ replaced by $\alpha^{(N)}_{n,m}$ and $n+m$ taken modulo $N$. Since $\alpha^{(N)}_{n,m}\to n\star m$ as $N\to\infty$, the scaled commutator converges to the Poisson bracket mode by mode. The second relation is kept for $n\in\Lambda_N$: the quantized Laplacian acts by $J_n\mapsto -\|n\|^2J_n$. Making these two replacements in the computation that leads to \eqref{eq:fourier} gives the Zeitlin equation:
\begin{equation}\label{eq_zeitlin_fourier}
\partial_t\omega_p =
\sum_{\substack{k\in\Lambda_N \\ k\neq p}}
\frac{1}{{\|k\|^2}}
\alpha^{(N)}_{p,k}\,\omega_k\,\omega_{p-k}, \qquad p\in\Lambda_N,
\end{equation}
where $p-k$ is reduced modulo $N$ to $\Lambda_N$.

\section{The Geometry of $\mathfrak{su}(N)$ in the Zeitlin Model}

In this section, we develop the Riemannian geometry of the finite-dimensional $\mathfrak{su}(N)$ arising from the Zeitlin model in order to compute its curvature explicitly. We first construct a convenient orthogonal basis of $\mathfrak{su}(N)$ and compute the corresponding adjoint and coadjoint operators, as well as the Levi-Civita connection. These are used to compute the sectional and Ricci curvatures, with the Ricci tensor diagonalizing in our chosen basis. The Ricci curvature of these finite-dimensional approximations will later be used to introduce a notion of Ricci curvature in $\mathrm{Diff}_{\mu, \mathrm{ex}}(T^2)$. 

Following the construction in \cite{Ze1991}, for each odd integer $N$, one has a collection of $N \times N$ complex matrices $J_{n}$, where $n = (n_1, n_2) \in \mathbb{Z}^2$, $n \neq (0, 0)$, satisfying
\begin{equation}\label{eq_J_definition}
J_nJ_m = \lambda^{(m\star n)/2} J_{n + m}, \qquad J_n^{\dagger} = J_{-n}.
\end{equation}
Here, $\lambda = e^{4\pi i/N}$ is a primitive $N$th root of unity, $m\star n := m_1n_2 - m_2n_1$, a skew-symmetric product in $\mathbb{Z}^2$, and $^{\dagger}$ is the conjugate transpose. Furthermore, the matrices $J_n$ exhibit some periodicity: if $n = m$ modulo $N$, then $J_n = J_m$, so we may think of the indices $n$ as elements of the group $\Z_N^2$ whenever this is convenient. When doing so, we shall choose symmetric representatives from the index set
\begin{equation}\label{Lamba_definition}
    \Lambda_N=\left\{(n_1,n_2)\in\Z^2\backslash \{0\}:-\frac{N-1}{2}\leq n_1,n_2,\leq \frac{N-1}{2}\right\}.
\end{equation}
It follows from \eqref{eq_J_definition} and the choice of $\lambda$ that
\begin{equation}\label{eq_lie_bracket}
[J_n, J_m] = -2i\sin\left( \frac{2\pi(n\star m)}{N}\right)\,J_{n+m}.
\end{equation}
To transform the $J_n$ into an $\mathfrak{su}(N)$ basis, define
\begin{equation*}
S_n = \frac{1}{2}\Big(J_n - J_{-n}\Big), \qquad C_n = \frac{1}{2i}\Big(J_n + J_{-n}\Big).
\end{equation*}
Note that the symmetries $C_{-n} = C_{n}$ and $S_{-n} = -S_{n}$ introduce some repetitions on the index set \eqref{Lamba_definition}. The \emph{positive} index set
\begin{equation*}
    \Lambda_N^+=\{(n_1,n_2)\in \Lambda_N: n_1 > 0\text{ or }n_1 = 0,\ n_2 > 0\}
\end{equation*}
eliminates these redundancies, so that the collection
\begin{equation*}
\big\{ C_n,\, S_n \, : \, (n_1, n_2) \in \Lambda_N^{+} \big\}
\end{equation*}
is a basis of $\mathfrak{su}(N)$. The quantization framework requires us to work with a scaled version of the standard structure constants arising from this basis. Thus, our structure constants are
\begin{equation}\label{eq_lie_bracket2}
\begin{aligned}
\ad_{S_n}S_m &= -\alpha^{(N)}_{n, m}\big(C_{n-m}+C_{n+m}\big),&\ad_{S_n}C_m &= -\alpha^{(N)}_{n, m}\big(S_{n-m}-S_{n+m}\big),
\\
\ad_{C_n}S_m&=-\alpha^{(N)}_{n,m}\big(-S_{n-m}-S_{n+m}\big),&\ad_{C_n}C_m &= -\alpha^{(N)}_{n, m}\big(C_{n-m}-C_{n+m}\big),
\end{aligned}
\end{equation}
where
\begin{equation*}
 \alpha^{(N)}_{n, m} := \frac N{2\pi}\sin\!\left(\frac{2\pi(n\star m)}{N}\right).
 \end{equation*} 
The \emph{Zeitlin metric}, introduced in \cite{Ze1991}, is
\begin{equation*}
\langle J_n, J_m \rangle = \delta_{n, m}\|n\|^2 = \delta_{n, m}\big(n_1^2 + n_2^2\big),
\end{equation*}
which in terms of our basis becomes 
\begin{equation}\label{eq:basisnorm}
\|S_n\|^2 = \|C_n\|^2 = \frac{1}{2}\|n\|^2, \quad \langle C_n, C_m \rangle = \langle S_n, S_m \rangle = \langle C_n, S_m \rangle = 0, \quad\text{for $n \neq m$}.
\end{equation}
Under $\sin(n\cdot x)\mapsto \tfrac{1}{2} S_n$ and $\cos(n\cdot x)\mapsto -\tfrac{1}{2} C_n$, the brackets~\eqref{eq_lie_bracket2} converge to the Poisson bracket, and~\eqref{eq:basisnorm} is the quantized analog of the $H^1$ metric on stream functions. It is also the $L^2$ metric on divergence-free vector fields, computed with normalized volume form so that $T^2$ has area $\tfrac{1}{4}$.

The corresponding coadjoint operator is:
\begin{align*}
\ad^*_{C_n} C_m &= \|m\|^2 \alpha^{(N)}_{n,m}\left(
\frac{C_{n-m}}{\|n-m\|^2} - \frac{C_{n+m}}{\|n+m\|^2}
\right),
&
\ad^*_{S_n} S_m &= \|m\|^2 \alpha^{(N)}_{n,m}\left(
\frac{C_{n-m}}{\|n-m\|^2} + \frac{C_{n+m}}{\|n+m\|^2}
\right), 
\\
\ad^*_{C_n} S_m &= -\|m\|^2 \alpha^{(N)}_{n,m}\left(
\frac{S_{n-m}}{\|n-m\|^2} + \frac{S_{n+m}}{\|n+m\|^2}
\right),
&
\ad^*_{S_n} C_m &= \|m\|^2 \alpha^{(N)}_{n,m}\left(
\frac{S_{n-m}}{\|n-m\|^2} - \frac{S_{n+m}}{\|n+m\|^2}
\right).
\end{align*}
Note that the structure constants vanish whenever $n=m$, and thus $\ad^*$ on any basis elements with $n=m$ is zero. 

\noindent Using the coadjoint operators $\ad^*$ computed above, we can also write the covariant derivatives: \begin{gather}
\begin{split}\label{eq_nablas}
\nabla_{C_n}C_m &= \alpha^{(N)}_{n, m}\big(\conn{m}{-n}C_{n-m} - \conn{m}{n}C_{n+m} \big), \\
\nabla_{S_n}S_m &= \alpha^{(N)}_{n, m}\big(\conn{m}{-n}C_{n-m} + \conn{m}{n}C_{n+m} \big), \\
\nabla_{C_n}S_m &= -\alpha^{(N)}_{n, m}\big(\conn{m}{-n}S_{n-m} + \conn{m}{n}S_{n+m} \big), \\
\nabla_{S_n}C_m &= \alpha^{(N)}_{n, m}\big(\conn{m}{-n}S_{n-m} - \conn{m}{n}S_{n+m} \big), 
\end{split}
\end{gather}
where
\begin{equation}\label{Gamma_def}
    \conn{m}{n} = \frac{1}{2}\left(1 + \frac{\|m\|^2 - \|n\|^2}{\|m+n\|^2} \right).
\end{equation}
 We choose not to absorb $\alpha^{(N)}_{n, m}$ into the connection constants $\conn{m}{n}$ because the explicit dependence on $\alpha^{(N)}_{n, m}$ is interesting from a quantization point of view. Next, the Riemann curvature tensor
\begin{equation*}
R(X, Y)Z = \nabla_{X}\nabla_{Y}Z - \nabla_{Y}\nabla_{X}Z - \nabla_{[X, Y]}Z
\end{equation*}
is computed in terms of these structure constants.

\begin{lemma}\label{lemma_riemann}
The Riemann curvature tensor of the Zeitlin metric on $SU(N)$, computed at the Lie algebra with respect to the basis $\{C_n, S_n\}$, has the following structure.
\begin{gather}\label{eq_big_R}
\begin{split}
    R(C_\ell,C_m)C_n&=R^1_{\ell,m,n}C_{\ell+m+n}-R^2_{\ell,m,n}C_{\ell+m-n}+R^3_{\ell,m,n}C_{\ell-m+n}-R^4_{\ell,m,n}C_{\ell-m-n}, \\
    R(C_\ell,C_m)S_n&=R^1_{\ell,m,n} S_{\ell+m+n}+R^2_{\ell,m,n} S_{\ell+m-n}+R^3_{\ell,m,n}S_{\ell-m+n}+R^4_{\ell,m,n}S_{\ell-m-n}, \\
    R(C_\ell,S_m)C_n&=R^1_{\ell,m,n} S_{\ell+m+n}-R^2_{\ell,m,n} S_{\ell+m-n}-R^3_{\ell,m,n}S_{\ell-m+n}+R^4_{\ell,m,n}S_{\ell-m-n}, \\
    R(C_\ell,S_m)S_n&=-R^1_{\ell,m,n}C_{\ell+m+n}-R^2_{\ell,m,n}C_{\ell+m-n}+R^3_{\ell,m,n}C_{\ell-m+n}+R^4_{\ell,m,n}C_{\ell-m-n},\\
    R(S_\ell,S_m)C_n&=-R^1_{\ell,m,n} C_{\ell+m+n}+R^2_{\ell,m,n} C_{\ell+m-n}+R^3_{\ell,m,n}C_{\ell-m+n}-R^4_{\ell,m,n}C_{\ell-m-n},\\
    R(S_\ell,S_m)S_n&=-R^1_{\ell,m,n} S_{\ell+m+n}-R^2_{\ell,m,n} S_{\ell+m-n}+R^3_{\ell,m,n}S_{\ell-m+n}+R^4_{\ell,m,n}S_{\ell-m-n},
\end{split}
\end{gather}
where the coefficients $R^{j}_{\ell, m, n}$ are given by
\begin{align*}
R^1_{\ell,m,n}&=\alpha^{(N)}_{m,n}\alpha^{(N)}_{\ell,m+n}\Gamma(n,m)\Gamma(m+n,\ell)-\alpha^{(N)}_{\ell,n}\alpha^{(N)}_{m,\ell+n}\Gamma(n,\ell)\Gamma(\ell+n,m)-\alpha^{(N)}_{\ell,m}\alpha^{(N)}_{\ell+m,n}\Gamma(n,\ell+m), \\
    R^2_{\ell,m,n}&=\alpha^{(N)}_{m,n}\alpha^{(N)}_{\ell,m-n}\Gamma(n,-m)\Gamma(m-n,\ell)-\alpha^{(N)}_{\ell,n}\alpha^{(N)}_{m,\ell-n}\Gamma(n,-\ell)\Gamma(\ell-n,m)-\alpha^{(N)}_{\ell,m}\alpha^{(N)}_{\ell+m,n}\Gamma(n,-\ell-m), \\
    R^3_{\ell,m,n}&=\alpha^{(N)}_{m,n}\alpha^{(N)}_{\ell,m-n}\Gamma(n,-m)\Gamma(m-n,-\ell)+\alpha^{(N)}_{\ell,n}\alpha^{(N)}_{m,\ell+n}\Gamma(n,\ell)\Gamma(\ell+n,-m)
    +\alpha^{(N)}_{\ell,m}\alpha^{(N)}_{\ell-m,n}\Gamma(n,\ell-m), \\
R^4_{\ell,m,n}&=\alpha^{(N)}_{m,n}\alpha^{(N)}_{\ell,m+n}\Gamma(n,m)\Gamma(m+n,-\ell)+\alpha^{(N)}_{\ell,n}\alpha^{(N)}_{m,\ell-n}\Gamma(n,-\ell)\Gamma(\ell-n,-m)+\alpha^{(N)}_{\ell,m}\alpha^{(N)}_{\ell-m,n}\Gamma(n,-\ell+m).
\end{align*}
In the formulas above, the terms for which $\alpha^{(N)}_{j, k} = 0$ should be omitted. 
\end{lemma} 
\begin{proof}
It follows immediately from our formulas \eqref{eq_nablas} for the covariant derivative that vanishing of $\alpha^{(N)}_{n, m}$ implies vanishing of the corresponding covariant derivatives, and consequently terms involving these in the Riemann tensor. We compute directly from the definition, using \eqref{eq_nablas} and the connection coefficients \eqref{Gamma_def}:
\begin{gather*}
\begin{split}
\nabla_{C_\ell}\nabla_{C_m}C_n&= \alpha^{(N)}_{m,n}\alpha^{(N)}_{\ell,m-n}\Gamma(n,-m)\Gamma(m-n,-\ell)C_{\ell-m+n}-\alpha^{(N)}_{m,n}\alpha^{(N)}_{\ell,m-n}\Gamma(n,-m)\Gamma(m-n,\ell)C_{\ell+m-n}\\
&\qquad-\alpha^{(N)}_{m,n}\alpha^{(N)}_{\ell,m+n}\Gamma(n,m)\Gamma(m+n,-\ell)C_{\ell-m-n}+\alpha^{(N)}_{m,n}\alpha^{(N)}_{\ell,m+n}\Gamma(n,m)\Gamma(m+n,\ell)C_{\ell+m+n}\\
-\nabla_{C_{m}}\nabla_{C_\ell}C_n&= -\alpha^{(N)}_{\ell,n}\alpha^{(N)}_{m,\ell-n}\Gamma(n,-\ell)\Gamma(\ell-n,-m)C_{\ell-m-n}+\alpha^{(N)}_{\ell,n}\alpha^{(N)}_{m,\ell-n}\Gamma(n,-\ell)\Gamma(\ell-n,m)C_{\ell+m-n}\\
&\qquad +\alpha^{(N)}_{\ell,n}\alpha^{(N)}_{m,\ell+n}\Gamma(n,\ell)\Gamma(\ell+n,-m)C_{\ell-m+n}-\alpha^{(N)}_{\ell,n}\alpha^{(N)}_{m,\ell+n}\Gamma(n,\ell)\Gamma(\ell+n,m)C_{\ell+m+n}\\
-\nabla_{[C_\ell,C_m]}C_n&= -\alpha^{(N)}_{\ell,m}\alpha^{(N)}_{\ell-m,n}\Gamma(n,-\ell+m)C_{\ell-m-n}+\alpha^{(N)}_{\ell,m}\alpha^{(N)}_{\ell-m,n}\Gamma(n,\ell-m)C_{\ell-m+n}\\
&\qquad +\alpha^{(N)}_{\ell,m}\alpha^{(N)}_{\ell+m,n}\Gamma(n,-\ell-m)C_{\ell+m-n}-\alpha^{(N)}_{\ell,m}\alpha^{(N)}_{\ell+m,n}\Gamma(n,\ell+m)C_{\ell+m+n}
\end{split}
\end{gather*}
Adding these together gives the first formula in \eqref{eq_big_R}. The remaining formulas in \eqref{eq_big_R} can be computed similarly.
\end{proof}

Note that not all possible combinations of inputs of the Riemann curvature tensor are given in Lemma \ref{lemma_riemann}, but the missing ones can be computed by symmetry. Next, we compute sectional curvatures
\begin{equation*}
K(X, Y) = \frac{\langle R(X, Y)Y, X \rangle}{\|X\|^2\|Y\|^2 - \langle X, Y \rangle^2}
\end{equation*}
in the orthogonal basis $\{ C_n, S_n \}$. 

\begin{lemma}\label{lemma_sectional}
The sectional curvatures of $2$-planes spanned by elements of the basis $\{C_n, S_n\}$ are:
\begin{gather}\label{eq:sectional}
\begin{split}
K^{(N)}(C_n, C_m) &= \frac{\big(\alpha^{(N)}_{n, m}\big)^2}{2\|n\|^2\|m\|^2}
\left(
4(\|n\|^2 + \|m\|^2) -3(\beta_+^2 + \beta_-^2)
+ \frac{(\|n\|^2-\|m\|^2)^2(\beta_+^2 + \beta_-^2)}{\beta_+^2\beta_-^2}
\right),
\end{split}
\end{gather}
where $n, m \in \Lambda_N^+$, with $n \neq m$ and $\beta_{\pm} = \|m \pm n\|$. Furthermore, we also have
\begin{equation*}
K^{(N)}(C_n, S_m) = K^{(N)}(S_n, S_m) = K^{(N)}(C_n, C_m),\quad \text{for $n \neq m$, and}\quad K^{(N)}(C_n, S_n) = 0.
\end{equation*}
\end{lemma}

\begin{proof} 
Starting from the formula in Lemma \ref{lemma_riemann} and taking inner product with $C_n$, note that only the first term survives, as $\langle C_{2m\pm n}, C_n \rangle = 0$ for all choices of $m$ and $n$. Thus, we compute the relevant connection coefficients. They are:
\begin{align*}
\Gamma(m, n) &= \frac{\beta_+^2 + \|m\|^2 - \|n\|^2}{2\beta_+^2},& \qquad \Gamma(m, -n) &= \frac{\beta_-^2 + \|m\|^2 - \|n\|^2}{2\beta_-^2},
\\
\Gamma(n + m, -m) &= \frac{\beta_+^2 + \|n\|^2 - \|m\|^2}{2\|n\|^2},& \qquad \Gamma(n - m, m) &= \frac{\beta_-^2 + \|n\|^2 - \|m\|^2}{2\|n\|^2}, \\ 
\Gamma(m, n)\Gamma(n+m, -m) &= \frac{\beta_+^4-(\|m\|^2-\|n\|^2)^2}{4\beta_+^2\|n\|^2},& \qquad 
\Gamma(m, -n)\Gamma(n-m, m) &= \frac{\beta_-^4-(\|m\|^2-\|n\|^2)^2}{4\beta_-^2\|n\|^2}, \\
\Gamma(m, n-m) &= \frac{\|n\|^2 + \|m\|^2 - \beta_-^2}{2\|n\|^2},& \qquad 
\Gamma(m, -n-m) &= \frac{\|n\|^2 + \|m\|^2 - \beta_+^2}{2\|n\|^2}. 
\end{align*}
Adding the last two, we get
\begin{equation*}
\Gamma(m, n-m) + \Gamma(m, -n-m) = \frac{2(\|n\|^2+\|m\|^2) -(\beta_+^2 + \beta_-^2)}{2\|n\|^2}.
\end{equation*}
Putting everything together,
\begin{gather*}
\begin{split}
&\Gamma(m, n)\Gamma(n+m, -m) + \Gamma(m, -n)\Gamma(n-m,m) - \Gamma(m, n-m) - \Gamma(m, -n-m) \\
&= \frac{\beta_+^4-(\|m\|^2-\|n\|^2)^2}{4\beta_+^2\|n\|^2}
+\frac{\beta_-^4-(\|m\|^2-\|n\|^2)^2}{4\beta_-^2\|n\|^2}
-\frac{2(\|n\|^2+\|m\|^2) -(\beta_+^2 + \beta_-^2)}{2\|n\|^2} \\
&= \frac{1}{4\|n\|^2}\left(
3(\beta_+^2+\beta_-^2)
-4(\|n\|^2+\|m\|^2)
-\frac{(\|n\|^2-\|m\|^2)^2(\beta_+^2+\beta_-^2)}{\beta_+^2\beta_-^2}
\right) 
\end{split}
\end{gather*}
after simplifying. Since the basis $\{C_n, S_n\}$ is orthogonal, $\langle R(C_n, C_m)C_m, C_n\rangle$ is equal to the above quantity scaled by $-\big(\alpha_{n, m}^{(N)}\big)^2\|C_n\|^2$, and the result for $K^{(N)}(C_n, C_m)$ follows. The identities for $K^{(N)}(S_n, S_m)$ and $K^{(N)}(C_n, S_m)$ are derived essentially in the same way, and $K^{(N)}(C_n, S_n) = 0$ because $\alpha_{n, n}^{(N)} = 0$.
\end{proof}

\begin{remark}\label{arnold_remark}
Fix $n, m \in \Lambda_N^+$ with $n \neq \pm m$. Then for $N$ sufficiently large, we have the identity
\begin{equation*}
\beta_+^2 + \beta_-^2 = 2(\|n\|^2 + \|m\|^2),
\end{equation*}
which can be used to simplify \eqref{eq:sectional}. In this case, letting $\theta_1$ denote the angle between the lattice vectors $m$ and $n$, we also have $\|n\|^2\|m\|^2\sin^2(\theta_1) = (m \star n)^2$. Then $K$ can be written as
\begin{gather*}
\begin{split}
K^{(N)}(C_n, C_m) &= -4(\alpha^{(N)}_{n, m})^2\sin^2(\theta_1)\left(
\frac{\big(\|m\|^2 + \|n\|^2\big)}{\|m-n\|^2\|m+n\|^2}\right).
\end{split}
\end{gather*}
In the large $N$ limit (taken along odd $N$), we have  $\alpha^{(N)}_{n, m} \rightarrow n \star m$. We thereby obtain
\begin{gather*}
\begin{split}
\lim\limits_{\substack{N \rightarrow \infty \\ \text{$N$ odd}}} K^{(N)}(C_n, C_m) &= -\sin^2(\theta_2)\sin^2(\theta_1)\big(\|m\|^2 + \|n\|^2\big),
\end{split}
\end{gather*}
where $\theta_2$ is the angle between $m+n$ and $m-n$, matching Arnold's computation \cite{arnold1966sur} up to scaling.
\end{remark}

\subsection{Ricci Curvature} A natural measure of the average sectional curvature experienced by a geodesic is captured by Ricci curvature, which makes it relevant to understanding the stability of fluid flows in the geometric framework of Arnold. We show that the Ricci tensor diagonalizes in the $\{C_n, S_n\}$ basis of $\mathfrak{su}(N)$, allowing its components to be computed explicitly.

\begin{theorem}
    The Ricci tensor diagonalizes in the $\{C_n, S_n\}$ basis of $\mathfrak {su}(N)$, with diagonal terms
\begin{equation}\label{eq_Ricci}
\mathrm{Ric}^{(N)}(C_n,C_n)=\mathrm{Ric}^{(N)}(S_n,S_n)=
\frac{1}{2}\sum_{\substack{m \in \Lambda_N \\ m \ne \pm n}}(\alpha^{(N)}_{n, m})^2\mathcal R(m,n),
\end{equation}
where, setting $\beta_{\pm} = \|m \pm n\|$ as before, 
\begin{equation*}\label{R_defn}
\mathcal R(m,n) := \frac{\|n\|^4-\big(\beta_+^2-\|m\|^2\big)^2}{\|m\|^2\beta_+^2}.
\end{equation*}
\end{theorem}
\begin{proof}
Recall that our index set $$\Lambda_N=\left\{(n_1,n_2)\in\Z^2\setminus \{0\}:-\frac{N-1}{2}\leq n_1,n_2\leq \frac{N-1}{2}\right\}$$
allows us to select representatives for our basis with indices in
$$\Lambda_N^+=\big\{(n_1,n_2)\in \Lambda_N: n_1 > 0\text{ or }n_1 = 0,\ n_2 > 0\big\}.$$
The diagonal terms of Ricci are then given by
\begin{align*}
    \mathrm{Ric}^{(N)}(C_n,C_n)&=\|C_n\|^2\sum_{\substack{m \in \Lambda_N^+\\m\neq n}} \Big(K^{(N)}(C_n,C_m)+K^{(N)}(C_n,S_m)\Big)\\
    &=\|C_n\|^2\sum_{\substack{m\in\Lambda_N\\m\neq \pm n}}K^{(N)}(C_n,C_m),
\end{align*}
where we use the fact that $K^{(N)}(C_n,S_m) = K^{(N)}(C_n,C_m)$ to re-index the sum from the positive index set $\Lambda_N^+$ to the full symmetric one $\Lambda_N$. Then~\eqref{eq:basisnorm} and ~\eqref{eq:sectional} allow us to write:
\begin{align*}
\mathrm{Ric}^{(N)}(C_n,C_n) &=
\frac{1}{4}\!\sum\limits_{\substack{m \in \Lambda_N \\ m \ne \pm n}}
\frac{(\alpha^{(N)}_{n, m})^2}{\|m\|^2}\left( 
4(\|n\|^2+\|m\|^2) - 3(\beta_+^2 + \beta_-^2)
+(\|n\|^2\!-\!\|m\|^2)^2\left(
\frac{1}{\beta_+^2} + \frac{1}{\beta_-^2}\right)
\right).
\end{align*}
Now note that the symmetry $m \mapsto -m$ leaves $(\alpha^{(N)}_{n, m})^2$ and $\|m\|^2$ unchanged, but swaps $\beta_+$ and $\beta_-$. Thus, since we are summing over the \emph{symmetric} index set $\Lambda_N$, we can effectively treat $\beta_+$ and $\beta_-$ as the same variable (by splitting the sum, applying the symmetry and recombining), which gives
\begin{align*}
\mathrm{Ric}^{(N)}(C_n,C_n) &=
\frac{1}{4} \sum\limits_{\substack{m \in \Lambda_N \\ m \ne \pm n}}
\frac{(\alpha^{(N)}_{n, m})^2}{\|m\|^2}\left( 
4(\|n\|^2+\|m\|^2) - 6\beta_+^2
+2\frac{(\|n\|^2-\|m\|^2)^2}{\beta_+^2}
\right) \\
&=
\frac{1}{2} \sum\limits_{\substack{m \in \Lambda_N \\ m \ne \pm n}}
\frac{(\alpha^{(N)}_{n, m})^2}{\|m\|^2}\left( 
2(\|n\|^2+\|m\|^2) - 3\beta_+^2
+\frac{(\|n\|^2-\|m\|^2)^2}{\beta_+^2}
\right) \\
&=
\frac{1}{2} \sum\limits_{\substack{m \in \Lambda_N \\ m \ne \pm n}}
\frac{(\alpha^{(N)}_{n, m})^2}{\|m\|^2\beta_+^2}\bigg( 
2\beta_+^2(\|n\|^2+\|m\|^2) - 3\beta_+^4 +
(\|n\|^2-\|m\|^2)^2
\bigg).
\end{align*}
Similarly, the symmetry $m \mapsto -m-n$ leaves $\|n\|^2$ and $(\alpha^{(N)}_{n, m})^2$ unchanged, but swaps $\|m\|^2$ with $\beta_+^2 = \|m+n\|^2$ (technically, this symmetry does not preserve the index set $\Lambda_N\backslash{}\{\pm n\}$, but the summand $m = n$ vanishes). We take the two sums, the original and the one after applying the symmetry $m \mapsto -m-n$, and average them. The relevant factor produced by this is:
\begin{align*}
\frac{1}{2}\bigg( 
\underbrace{2\beta_+^2(\|n\|^2+\|m\|^2) - 3\beta_+^4 +
(\|n\|^2-\|m\|^2)^2}_{\text{original}}
+
\underbrace{2\|m\|^2(\|n\|^2+\beta_+^2) - 3\|m\|^4 +
(\|n\|^2-\beta_+^2)^2}_{\text{swapped: $\|m\| \leftrightarrow \beta_+$}}
\bigg)
\end{align*}
which upon simplification gives the formula claimed \eqref{eq_Ricci}. This takes care of the diagonal terms ${\mathrm{Ric}^{(N)}(C_n,C_n)}$, and the (also diagonal) $\mathrm{Ric}^{(N)}(S_n,S_n)$ terms are computed similarly.

For off-diagonal terms involving $C_n$ and $S_m$ for any $n,m,$ we have
\begin{equation}\label{eq_ricci_mixed_1}
\mathrm{Ric}^{(N)}(C_n, S_m) = \sum\limits_{\substack{\ell \in \Lambda_N^+ \\ \ell \ne m, n}}
\frac{2}{\|\ell\|^2}\bigg(\big\langle
R(C_{\ell}, C_n)S_m, C_{\ell}
\big\rangle
+
\big\langle
R(S_{\ell}, C_n)S_m, S_{\ell}
\big\rangle\bigg).
\end{equation}
We claim that both summands in \eqref{eq_ricci_mixed_1} vanish. For the first one, the covariant derivative formulas computed in \eqref{eq_nablas} show that
\begin{equation*} 
R(C_\ell, C_n)S_m=\nabla_{C_\ell}\nabla_{C_n}S_m - \nabla_{C_n}\nabla_{C_\ell}S_m - \nabla_{[C_\ell, C_n]}S_m
\end{equation*} 
produces a linear combination of the $S$ basis elements, and since $\langle S_{\ell}, C_{\ell'} \rangle = 0$ for any choice of indices $\ell, \ell'$ by construction,
\begin{equation*} 
\big\langle R(C_\ell, C_n)S_m, C_\ell\big\rangle = 0.
\end{equation*}
Likewise, the expression 
\begin{equation*} 
R(S_\ell, C_n)S_m = \nabla_{S_\ell}\nabla_{C_n}S_m - \nabla_{C_n}\nabla_{S_\ell}S_m - \nabla_{[S_\ell, C_n]}S_m
\end{equation*}
produces a linear combination of the $C$ basis elements, so $\langle R(S_\ell, C_n)S_m, S_\ell\rangle = 0$. This takes care of the ``mixed'' (involving $C$ and $S$ basis elements) Ricci terms.

The remaining off-diagonal terms are $(C_n, C_m)$ and $(S_n, S_m)$ for $n \neq \pm m$. Once again the computations are very similar between $C$ or $S$, so we work out only the $C$ case explicitly, which is
\begin{equation*}
\mathrm{Ric}^{(N)}(C_n, C_m) =  \sum\limits_{\substack{\ell \in \Lambda_N^+ \\ \ell \ne m, n}}
\frac{2}{\|\ell\|^2}\bigg(\big\langle
R(C_{\ell}, C_n)C_m, C_{\ell}
\big\rangle
+
\big\langle
R(S_{\ell}, C_n)C_m, S_{\ell}
\big\rangle\bigg).
\end{equation*}
In the formulas below, connection contributions with vanishing structure constants are taken to be zero by convention, without evaluating the associated $\Gamma$ expressions, and the terms whose output index is zero modulo $N$ are omitted.

The expression for $R(C_{\ell}, C_n)C_m$ involves only the four basis vectors $C_{\ell\pm m\pm n}$, and that of $R(S_{\ell}, C_n)C_m$ involves only the basis elements $S_{\ell\pm m\pm n}$. Specifically, they read
\begin{align*}
R(C_\ell,C_n)C_m&=R^1_{\ell,n,m}C_{\ell+n+m}-R^2_{\ell,n,m}C_{\ell+n-m}+R^3_{\ell,n,m}C_{\ell-n+m}-R^4_{\ell,n,m}C_{\ell-n-m},\\
R(S_\ell,C_n)C_m&=R^1_{\ell,n,m}S_{\ell+n+m}-R^2_{\ell,n,m}S_{\ell+n-m}+R^3_{\ell,n,m}S_{\ell-n+m}-R^4_{\ell,n,m}S_{\ell-n-m},
\end{align*}
where the coordinate functions are given by \eqref{eq_big_R}, 
with $\Gamma$ defined in~\eqref{Gamma_def}. 

The vectors $R(C_\ell,C_n)C_m$ and $R(S_\ell,C_n)C_m$ thereby survive the inner product with $C_\ell$ and $S_\ell$, respectively, only when $\ell\pm m\pm n\equiv \pm \ell\pmod N,$ which forces either $m = \pm n$ or $2\ell = (\pm m\pm n) \pmod N$. The former case constitutes the diagonal terms of Ricci. In the latter case, non-vanishing terms involve $S_{-\ell}=-S_\ell$ and $C_{-\ell}=C_\ell.$ But for any such term, one has 
$$\big\langle R(C_{\ell}, C_n)C_m, C_{\ell}\big\rangle+\big\langle
R(S_{\ell}, C_n)C_m, S_{\ell}\big\rangle=(-1)^{j+1}(R^j_{m\ell n}\|C_\ell\|^2-R^j_{m\ell n}\|S_\ell\|^2)=0$$
for some $j\in\{1,2,3,4\}.$ Therefore, $\mathrm{Ric}(C_n, C_m)=0$ for $m \neq \pm n.$ An analogous argument shows that $\mathrm{Ric}(S_n, S_m)=0$ when $m \neq \pm n.$ 
\end{proof}

\begin{remark}
Since $\|C_n\|^2 = \|S_n\|^2 = \tfrac{1}{2}\|n\|^2$, this proves Theorem \ref{eq_intro_thm_1}.
\end{remark}

The diagonalization of the Ricci tensor allows the curvature in the direction of each basis element to be computed independently as a finite sum involving the structure constants of the Zeitlin Lie algebra. In the next section, we use this expression to study the behavior of the curvature as the dimension of our approximation increases. 

\section{Convergence of Averaged Ricci Curvature}

In this section we prove Theorem~\ref{eq_main_thm_intro}, the convergence of averaged Ricci curvature in the large $N$ limit. Afterwards, we discuss the asymptotic expansion \eqref{torus_asymp} mentioned in the introduction. 

\begin{proof}[Proof of Theorem \ref{eq_main_thm_intro}] Fix $n\in\mathbb{Z}^2\setminus\{0\}$, and let $N$ range over odd integers large enough that $n\in\Lambda_N$. By Theorem~\ref{eq_intro_thm_1} and~\eqref{eq:basisnorm}, $r_n(N) = 2\,\mathrm{Ric}^{(N)}(C_n,C_n)/\|n\|^2$, so that
\begin{equation*}
\widetilde r_n(N) = \frac{2N^2}{(N^2-1)\,\|n\|^2}\frac{1}{N^2}\,\mathrm{Ric}^{(N)}(C_n,C_n).
\end{equation*}
Since $N^2/(N^2-1)\to 1$, Theorem~\ref{eq_main_thm_intro} is equivalent to
\begin{equation*}\label{eq_lim_4}
\lim_{N\to\infty}\frac{1}{N^2}\,\mathrm{Ric}^{(N)}(C_n,C_n) = -\frac{1}{2\pi^2}\int_Q \frac{\sin^2\big(2\pi(n\star x)\big)\,(n\cdot x)^2}{\|x\|^4}\,dx, \quad Q = \big[-\tfrac{1}{2},\tfrac{1}{2}\big]^2.
\end{equation*}
Our starting point is the formula
\begin{equation*}\label{eq_lim_5}
\frac{1}{N^2}\,\mathrm{Ric}^{(N)}(C_n,C_n) = \frac{1}{2N^2}\sum_{\substack{m\in\Lambda_N\\ m\neq\pm n}} \big(\alpha^{(N)}_{n,m}\big)^2\,\mathcal{R}(m,n),
\end{equation*}
where, as before, $\mathcal R(m,n)$ and the scaled structure constants $\alpha^{(N)}_{n, m}$ are given by
\begin{align*}
\mathcal{R}(m,n) &= \frac{\|n\|^4-\big(\beta_+^2-\|m\|^2\big)^2}{\|m\|^2\beta_+^2}, \qquad \beta_+ = \|m+n\|, \\
\alpha_{n,m}^{(N)} &= \frac{N}{2\pi}
\sin\left(\frac{2\pi(n\star m)}{N}\right).
\end{align*}
Next, partition the square $Q = \big[-\tfrac{1}{2}, \tfrac{1}{2}\big]^2$ into $N^2$ half-open (so they are disjoint) cells uniformly, each having area $1/N^2$. Any $x \in Q$ belongs to exactly one of these cells, except for a measure zero set on the boundary. For $x \in \mathrm{int}(Q)$, let $x_N$ be the \emph{center} of the unique cell containing $x$, which is
\begin{equation*}
x_N = \frac{1}{N}\Big(\floor{Nx_1 + \tfrac{1}{2}}, ~\floor{Nx_2 + \tfrac{1}{2}}\Big).
\end{equation*}
We define a sequence of functions $f_N : Q \rightarrow \mathbb{R}$ by setting
\begin{gather*}\label{eq_def_fN}
\begin{split}
f_N(x) := -\frac{1}{2}\big(\alpha^{(N)}_{n, m}\big)^2 \mathcal{R}\left(m, n\right), \qquad m = Nx_N,
\end{split}
\end{gather*}
for $x \in \mathrm{int}(Q)$ and $m \neq 0, \pm n$, and $f_N(x) := 0$ otherwise.

For every $N$, the function $f_N$ is the sum of $N^2$ step functions, since it is clearly constant on each of the $N^2$ cells of $Q$. Therefore, the integral of $f_N$ over the square $Q$ is
\begin{gather*}\label{eq_fN_integral}
\begin{split} 
\int\limits_Q f_N(x)\,dx
&=
\sum\limits_{m \in \tilde\Lambda_N} \frac{1}{N^2}\, f_N\left(\frac{m}{N}\right) \\ 
&= -\frac{1}{2N^2} \sum\limits_{\substack{
m \in \Lambda_N \\ m \neq \pm n
}} (\alpha^{(N)}_{n, m})^2\mathcal{R}(m, n) \\
&= -\frac{1}{N^2} \mathrm{Ric}^{(N)}(C_n, C_n)
\end{split} 
\end{gather*}
where $\tilde \Lambda_N=\Lambda_N\cup\{0\}.$ Note that the cells indexed by $m = \pm n$ and $m = 0$ contribute zero by the definition of $f_N(x)$.

As $N \rightarrow \infty$, for $x \in \mathrm{int}(Q)\backslash{}\{0\}$ we get the pointwise limit
\begin{equation}\label{eq_def_f}
\lim\limits_{N \rightarrow \infty} f_N(x) = f(x),\quad\text{where}\quad f(x) := \frac{1}{2\pi^2}\sin^2\big(2\pi(n \star x)\big)\frac{(n\cdot x)^2}{\|x\|^4},
\end{equation}
where $f(0)=0$. To prove \eqref{eq_def_f}, fix $x \in \mathrm{int}(Q)\backslash{}\{0\}$. Then $x_N \rightarrow x$ and $n/N \rightarrow 0$, so for sufficiently large $N$, the index $m+n$ need not be reduced modulo $N$, and we can use the identity
\begin{equation*}
\|n\|^4 - (\beta_+^2 - \|m\|^2)^2 = \|n\|^4 - (2n\cdot m + \|n\|^2)^2 = -4(n\cdot m)(n\cdot m + \|n\|^2),
\end{equation*}
which simplifies the numerator of $\mathcal{R}(m, n)$. Plugging this into $f_N(x)$ gives
\begin{equation*}
f_N(x) = \frac{\sin^2\big(2\pi(n \star x_N)\big)}{2\pi^2}\frac{(n \cdot x_N)^2 + \frac{\|n\|^2}{N}(n\cdot x_N)}{\|x_N\|^2\left\|x_N + \frac{n}{N}\right\|^2},
\end{equation*}
and taking the limit yields \eqref{eq_def_f}. This proves pointwise convergence almost everywhere, since the boundary of $Q$ (together with the origin) has measure zero.

Thus, if we can show that the sequence $f_N$ is dominated by an $L^1$ function $g$ on $Q$, it follows by the dominated convergence theorem that
\begin{equation*}
\lim\limits_{N \rightarrow \infty} \frac{1}{N^2} \mathrm{Ric}^{(N)}(C_n, C_n) = -\lim\limits_{N \rightarrow \infty} \int\limits_Q f_N(x)\,dx = -\int\limits_Q f(x)\,dx.
\end{equation*}
We claim that
\begin{equation}\label{eq_fN_bound}
|f_N(x)| \leq 2\|n\|^4.
\end{equation}
This is immediate on the boundary $\partial Q$ and on the cells indexed by $m \in \{0, \pm n\}$, since $f_N = 0$ there by definition. We prove \eqref{eq_fN_bound} for the other, nontrivial cases as follows.

The standard inequality $|\sin(x)| \leq |x|$ gives
\begin{equation*} 
\big(\alpha^{(N)}_{n, m}\big)^2 \leq \|n\|^2\|m\|^2.
\end{equation*}
Note that since $\alpha^{(N)}_{n, m} = \alpha^{(N)}_{n, m+n}$, we also have
\begin{equation*} 
\big(\alpha^{(N)}_{n, m}\big)^2 \leq \|n\|^2\|m+n\|^2,
\end{equation*}
so putting these together gives
\begin{equation*} 
\big(\alpha^{(N)}_{n, m}\big)^2 \leq \|n\|^2\, \mathrm{min}\big\{\|m\|^2,\, \|m+n\|^2\big\}.
\end{equation*}
Here we point out that the triangle inequality $\|m+n\| \leq \|m\| + \|n\|$ still holds even as the multi-indices $m$ and $n$ are reduced modulo $N$, since reduction can only make $\|m+n\|$ smaller (with our symmetric choice of representatives $\Lambda_N$). From this we have the standard inequalities
\begin{equation*}
\big| \|m + n\| - \|m\| \big| \leq \| n \| \leq \|m + n\| + \|m\|,
\end{equation*}
which give the bound
\begin{align*}
\big| \|n\|^4 - (\|m+n\|^2 - \|m\|^2)^2 \big|
&\leq
\mathrm{max}\big\{
\|n\|^4, \, (\|m+n\|^2 - \|m\|^2)^2
\big\} \\
&\leq \|n\|^2\,\big(\|m+n\| + \|m\|\big)^2 \\
&\leq 4\|n\|^2\,\mathrm{max}\big\{\|m+n\|^2, \, \|m\|^2\big\}.
\end{align*}
From the above inequalities it follows that
\begin{align*}
|f_N(x)| &= \frac{(\alpha^{(N)}_{n, m})^2}{2\|m+n\|^2\|m\|^2}\,
\Big|\|n\|^4 - \big(\|m+n\|^2 - \|m\|^2\big)^2
\Big| \\
&\leq \frac{\|n\|^2\, \mathrm{min}\big\{\|m\|^2,\, \|m+n\|^2\big\}}{2\|m+n\|^2\|m\|^2}\,
4\|n\|^2\,\mathrm{max}\big\{\|m+n\|^2, \, \|m\|^2\big\} \\
&=2\|n\|^4 \bigg(\frac{\mathrm{min}\big\{\|m\|^2,\, \|m+n\|^2\big\}\, \mathrm{max}\big\{\|m+n\|^2, \, \|m\|^2\big\}}{\|m+n\|^2\|m\|^2}\bigg)
\\
&= 2\|n\|^4.
\end{align*}
This is the claimed bound, so dominated convergence gives
\begin{gather}\label{eq_integral_rho}
\begin{split}
\lim\limits_{\substack{N \rightarrow \infty \\ \text{$N$ odd}}} \frac{1}{N^2} \mathrm{Ric}^{(N)}(C_n, C_n)
&= -\lim\limits_{\substack{N \rightarrow \infty \\ \text{$N$ odd}}}
\int\limits_Q f_N(x)\,dx = -\int\limits_Q f(x)\,dx \\
&= -\frac{1}{2\pi^2}\int_Q \frac{\sin^2\big(2\pi(n \star x)\big)(n \cdot x)^2}{\|x\|^4}\,dx,
\end{split}
\end{gather}
which is strictly negative.
\end{proof}

\begin{proposition}\label{prop_asymp}
The averaged Ricci curvature $\rho_{T^2}(n)$ of $\mathrm{Diff}_{\mu, \mathrm{ex}}(T^2)$, given by the integral in \eqref{eq_limit_intro}, admits the asymptotic expansion
\begin{equation}\label{eq_asymp}
2\pi\rho_{T^2}(n) = -\log\|n\| - \gamma + \frac{1}{2} + \frac{2G}{\pi} - \log(2\pi) + o(1)
\end{equation}
as $\|n\| \rightarrow \infty$, where $\gamma$ is the Euler--Mascheroni constant and $G$ is Catalan's constant. 
\end{proposition}

\begin{proof}
We introduce polar coordinates $x = r(\cos\varphi, \sin\varphi)$, so that $dx = r\,dr\,d\varphi$ and
\begin{equation*} 
Q = \left[-\frac{1}{2},\,\frac{1}{2}\right]^2 = \big\{ (r, \varphi) \, : \, 0 \leq r \leq b(\varphi) \big\},
\end{equation*}
where $b(\varphi)$ is the distance from the origin to $\partial Q$ in the direction $\varphi$. Thus,
\begin{equation}\label{eq_b_def}
\frac{1}{b(\varphi)} = 2\max\{|\cos\varphi|,|\sin\varphi|\},
\end{equation} 
and $\tfrac{1}{2}\leq b(\varphi) \leq \tfrac{\sqrt2}{2}$. Let $\theta$ be the counterclockwise angle from $n$ to $x$, so that
\begin{equation}\label{eq_theta_def}
n\cdot x = \|n\|\,r\cos\theta, \quad n\star x = \|n\|\,r\sin\theta, \quad\text{and}\quad \cos\theta = \frac{n_1\cos\varphi + n_2\sin\varphi}{\|n\|}.
\end{equation}
Note that $\theta$ depends on $\varphi$. Then, we can write $\rho_{T^2}(n)$ as
\begin{equation*}
\rho_{T^2}(n) = -\frac{1}{\pi^2}\int\limits_0^{2\pi}\cos^2\theta\int\limits_0^{b(\varphi)}\frac{\sin^2\big(2\pi\|n\|\,r\sin\theta\big)}{r}\,dr\,d\varphi.
\end{equation*}
We evaluate the radial integral first. The identity $\sin^2(x) = \tfrac{1}{2}(1 - \cos(2x))$ gives, for any $a > 0$,
\begin{equation*}
\int\limits_0^{b}\frac{\sin^2(ar)}{r}\,dr = \frac{1}{2}\big(\gamma + \log(2ab) - \mathrm{Ci}(2ab)\big), \quad\text{where}\quad \mathrm{Ci}(z) := -\int\limits_z^\infty \frac{\cos t}{t}\,dt
\end{equation*}
is the cosine integral (see \cite[\S5.2]{abramowitz1964handbook}). Taking $a = 2\pi\|n\|\,|\sin\theta|$ and writing $z(\varphi) := 4\pi\|n\|\,b(\varphi)\,|\sin\theta|$, 
\begin{equation*}
\rho_{T^2}(n) = -\frac{1}{2\pi^2}\int\limits_0^{2\pi}\cos^2\theta\,\Big(\gamma + \log(4\pi\|n\|) + \log b(\varphi) + \log|\sin\theta| - \operatorname{Ci}\big(z(\varphi)\big)\Big)\,d\varphi.
\end{equation*}
We integrate the five terms above separately. The first two are constant, and we have $\int_0^{2\pi}\cos^2\theta\,d\varphi = \pi$ from the definition \eqref{eq_theta_def} of $\theta$. The third one is
\begin{equation*}
\int\limits_0^{2\pi}\cos^2\theta\,\log b(\varphi)\,d\varphi = -\frac{1}{\|n\|^2}\int\limits_0^{2\pi}\big(n_1\cos\varphi + n_2\sin\varphi\big)^2\Big(\log 2 + \log\max\big\{|\cos\varphi|,|\sin\varphi|\big\}\Big)\,d\varphi,
\end{equation*}
by \eqref{eq_b_def} and \eqref{eq_theta_def}. Splitting $[0, 2\pi]$ into the eight intervals of length $\pi/4$ on which the max is either $|\cos\varphi|$ or $|\sin\varphi|$, this reduces to integrals of $\log|\cos\varphi|$ and $\log|\sin\varphi|$ multiplied by $\cos^2\varphi$, $\sin\varphi\cos\varphi$ and $\sin^2\varphi$. Most of these are elementary, except for
\begin{equation*}
\int\limits_0^{\pi/4}\log\cos\varphi\,d\varphi = -\frac{\pi}{4}\log 2 + \frac{G}{2},
\end{equation*}
which can be found in standard tables (see \cite{gradshteyn2014table}).

For the fourth term, $\theta$ differs from $\varphi$ by a constant, so we may integrate in $\theta$ over a full period:
\begin{equation*}
\int\limits_0^{2\pi}\cos^2\theta\,\log|\sin\theta|\,d\varphi = \int\limits_0^{2\pi}\cos^2\theta\,\log|\sin\theta|\,d\theta = -\pi\log 2 - \frac{\pi}{2},
\end{equation*}
which is again standard: it follows from Euler's integral $\int_0^{\pi}\log\sin\theta\,d\theta = -\pi\log 2$ after writing $\cos^2\theta = \frac12(1+\cos2\theta)$ and integrating the $\cos 2\theta$ term by parts.

Finally, $\mathrm{Ci}$ is bounded on $[1,\infty)$, goes to $0$ at infinity, and satisfies $|\mathrm{Ci}(z)|\lesssim (1+|\log z|)$ on $(0,1]$. Since $b \geq \tfrac{1}{2}$, we have $z(\varphi) \geq |\sin\theta|$ for $\|n\|\ge1$, so the last term is dominated by the integrable function $1+|\log|\sin\theta||$ (up to constants), and goes to zero wherever $\sin\theta\neq0$ as $\|n\| \to \infty$. By dominated convergence, its contribution is $o(1)$. Collecting terms,
\begin{equation*}
\rho_{T^2}(n) = -\frac{1}{2\pi^2}\Big(\pi\gamma + \pi\log(4\pi\|n\|) - 2G - \pi\log 2 - \frac{\pi}{2}\Big) + o(1),
\end{equation*}
and multiplying by $2\pi$ gives~\eqref{eq_asymp}.
\end{proof}

\section{Conjugate Points and the Diameter of $SU(N)$}

In this section, we derive a diameter estimate for $SU(N)$, equipped with the Zeitlin metric and the scaled bracket \eqref{eq_lie_bracket2}. It is known \cite{lebrigant2025geodesics} that every geodesic of $SU(N)$, under any left or right-invariant metric, develops a conjugate point. Since geodesics cannot minimize distance beyond their first conjugate point, the core problem is to obtain a uniform upper bound on the first conjugate time for all geodesics, depending only on $N$.

The results of \cite{lebrigant2025geodesics} show that, for a geodesic $\gamma(t)$ starting at the identity with initial velocity\footnote{Non-zero velocity is a standing assumption for all geodesics throughout the paper, though $u_0$ need not be normalized.} $u_0 \neq 0$, the first conjugate time $t = \tau$ satisfies the bounds
\begin{equation}\label{conjugatetime}
\frac{2\pi\lambda_{min}(\Lambda)}{\|\ad_{\Lambda u_0}\|_{op}}\le\tau\le\frac{2\pi\lambda_{max}(\Lambda)}{\|\ad_{\Lambda u_0}\|_{op}},
\end{equation}
where the operator norms are computed with respect to the bi-invariant metric $\kappa$. Here, $\lambda_{min}$ and $\lambda_{max}$ are the smallest and largest eigenvalues of the operator 
\begin{equation*}
\Lambda : \mathfrak{su}(N) \rightarrow \mathfrak{su}(N)
\end{equation*}
relating the bi-invariant metric $\kappa$ to the Zeitlin metric $\langle\cdot,\cdot\rangle$. Recall that $\kappa(u,v) = -\tr(uv)$ is a bi-invariant metric on $SU(N)$. It is essentially the negative of the Killing form. Since any left- or right-invariant metric on $SU(N)$ is then related to $\kappa$ by some symmetric, positive-definite linear operator, we take $\Lambda$ to be the operator so that
\begin{equation}\label{determinesLambda}
    \langle u,v \rangle = \kappa(u, \Lambda v) =  -\tr(u\Lambda v),\qquad\Lambda:\mathfrak{su}(N)\to\mathfrak{su}(N).
\end{equation}

The diameter estimate is now obtained as follows. Every $g \in SU(N)$ can be joined to the identity $e$ by a minimizing geodesic, which we may take to have unit speed so that $\|u_0\| = 1$. A geodesic cannot minimize distance past its first conjugate point, so $\mathrm{dist}(g, e) \leq \tau(u_0)$. Moreover, right-invariance gives $\mathrm{dist}(g, h) = \mathrm{dist}(e, hg^{-1})$, and hence $\mathrm{diam}(SU(N)) = \sup_g \mathrm{dist}(e, g)$. Combining this with \eqref{conjugatetime} yields
\begin{equation}\label{diameter}
\operatorname{diam}\big(SU(N)\big)\le \frac{2\pi\lambda_{max}(\Lambda)}{\inf_{\|u_0\|=1}\|\ad_{\Lambda u_0}\|_{op}}.
\end{equation}
We begin by computing the operator norms in \eqref{conjugatetime} first in the case where $u_0$ is one of the basis vectors $C_n$ or $S_n$, and then use it to estimate the operator norm for general $u_0$.

\begin{lemma}
Let $n \in \Lambda_N^+$ and $d = \gcd(n_1, n_2, N)$. The operator norms of $\ad_{\Lambda C_n}$ and $\ad_{\Lambda S_n}$, computed with respect to the bi-invariant metric $\kappa(\cdot,\cdot)$ in $\mathfrak{su}(N)$, are
\begin{equation*} 
\|\ad_{\Lambda C_n}\|_{op} = \frac{\| n \|^2}{\pi}\cos^2\left(\frac{\pi d}{2N}\right),
\quad\text{and}\quad
\|\ad_{\Lambda S_n}\|_{op}=\frac{\|n\|^2}{\pi}\cos\Big(\frac{\pi d}{2N}\Big).
\end{equation*}
\end{lemma}

\begin{proof}
To start, we need the eigenvalues of the operator $\Lambda:\mathfrak{su}(N)\to \mathfrak{su}(N)$, which may be found using the relation~\eqref{determinesLambda} and the fact that 
\begin{equation}\label{Zeitlin_norms}
\|C_n\|^2=\|S_n\|^2=\frac{1}{2}\|n\|^2.
\end{equation}
Specifically, the multiplication rule~\eqref{eq_J_definition} for the Zeitlin basis $J_n$ of the complexification $\mathfrak{su}(N)_{\mathbb{C}}\cong\mathfrak{sl}(N,\mathbb{C})$ gives 
\begin{equation}\label{biinvariant_norms}
-\tr(S_nS_m)=-\tr(C_nC_m)=\frac{N}{2}\delta_{n,m}\quad\text{and}\quad -\tr(C_nS_m)=0.
\end{equation}
The basis $\{C_n,S_m\}$ of $\mathfrak{su}(N)$ is therefore orthogonal with respect to both the Zeitlin metric and the bi-invariant metric $\kappa$. Consequently, $\{C_n,S_m\}$ form an eigenbasis of $\Lambda$, so there are constants $\lambda_n$ and $\mu_n$ such that $\Lambda C_n=\lambda_n C_n$ and $\Lambda S_n=\mu_n S_n.$ Substituting these into the relation~\eqref{determinesLambda} and using~\eqref{Zeitlin_norms} and~\eqref{biinvariant_norms} shows that the eigenvalues are $\lambda_n=\mu_n=\|n\|^2/N,$ and hence
\begin{equation*}
    \Lambda C_n=\frac{\|n\|^2}{N}C_n,\qquad \Lambda S_n=\frac{\|n\|^2}{N}S_n.
\end{equation*}
It follows from these expressions and the range of permissible indices that
\begin{equation}\label{Lambdaeigenvalues}
    \lambda_{min}(\Lambda)=\frac{1}{N},\qquad\lambda_{max}(\Lambda)=\frac{(N-1)^2}{2N}.
\end{equation}
We will now calculate the desired operator norm, which must be done using the bi-invariant metric $\kappa(X,Y)=-\tr(XY)$. Since $L=\ad_{u_0}$ is skew-adjoint with respect to $\kappa$, the operator norm is given by 
\[\|\ad_{u_0}\|_{op}=\sqrt{\lambda_{max}(-L^2)}.\]
Using~\eqref{eq_lie_bracket2}, we have 
\begin{align*}
    -\ad_{C_n}\ad_{C_n}C_m
    &=(\alpha_{n,m}^{(N)})^2(2C_m-C_{m-2n}-C_{m+2n}),\\
    -\ad_{C_n}\ad_{C_n}S_{m}
    &=(\alpha_{n,m}^{(N)})^2(2S_{m}-S_{m-2n}-S_{m+2n}),
\end{align*}
from which it is clear that the matrix representation of $-\ad_{C_n}^2$ is naturally block-diagonal with two identical blocks: one for the $C_m$ basis elements, and one for the $S_m$ basis elements. Our goal is to construct this matrix in a way that will readily yield its eigenvalues. Since there are two identical blocks, we can focus solely on ordering the $C_m$ elements.

Fix $n=(n_1,n_2)\in\Z_N\times\Z_N$. We order the $C_m$ basis by first partitioning the indices $m\in\Z_N\times\Z_N$ into cosets of $\langle n \rangle$. Let $d=\gcd(n_1,n_2,N).$ Then $\langle n\rangle \subset \Z_N\times \Z_N$ has $Nd$ cosets of the form $\{m,m+n,m+2n,...,m+(N/d-1)n\}$, each of which has cardinality $N/d$.

We will not need to include every coset due to the property $C_m=C_{-m}$. In particular, nontrivial cosets come in pairs $\{m+\langle n\rangle,-m+\langle n\rangle\}.$ Indeed, $-m\notin m+\langle n \rangle$, else $2m\equiv kn\mod N$ for some $k<N/d$, where the invertibility of 2 in $\Z_N$ for odd $N$ forces $m\in\langle n\rangle,$ contradicting the nontriviality of $m+\langle n\rangle$. We then need only include one coset from each pair. This leaves a total of $\frac12(Nd-1)$ nontrivial cosets, each of which contains $N/d$ indices. Thus, we include $\frac12(Nd-1)(N/d)$ indices from the nontrivial cosets.

Moreover, the trivial coset $\langle n \rangle$ contributes $\frac12(N/d-1)$ indices, as it contains both the zero index and the additive inverse of each of its elements, which must be removed.

We have therefore accounted for the entire $C_m$ basis, since the above considerations give a total of  
\[\frac12 \left(\frac Nd-1\right)+\frac12\left(Nd-1\right)\left(\frac Nd\right)=\frac12(N^2-1)\]
indices included in this construction.

Grouping the indices of our $C_m$ basis elements by cosets of $\langle n\rangle$ yields a block-circulant matrix representation for the $C_m$ portion of $-\ad_{C_n}^2$, which takes the form
\begin{equation}\label{ad^2matrix}
\begin{bmatrix}
~\big[0\big] & & & & \\
&(\alpha^{(N)}_{n,m_1})^2M\\
&&\ddots&\\
&&&&(\alpha^{(N)}_{n,m_j})^2M
\end{bmatrix},
\end{equation}
where $j = \frac{1}{2}(Nd-1)$, the numbers $m_1,...,m_j$ are representatives of the $\frac{1}{2}(Nd-1)$ nontrivial cosets of $\langle n\rangle,$ and $M$ is the $\frac Nd\times \frac Nd$ circulant matrix
$$M=\begin{bmatrix}
    2&0&-1&0&\cdots&-1&0\\
    0&2&0&-1&&&-1 \\
    -1&0&2&0&\ddots && \vdots\\
    0&-1&0&
    2&&\ddots&\vdots\\
    \vdots&&\ddots&&\ddots&&-1\\
    -1&&&\ddots&&2&0\\
    0&-1&\dots&\dots&-1&0&2
\end{bmatrix}.$$
Here, the $k^{\text{th}}$ column of the $i^{\text{th}}$ circulant block is associated with the basis element $C_{m_i+(k-1)n}$. Note that the trivial coset $\langle n\rangle$ forms the zero block, since the structure constants $\alpha^{(N)}_{n, m}$ vanish. The eigenvalues $\mu_k$ of the matrix~\eqref{ad^2matrix} are, up to the factors $(\alpha_{n, m_i}^{(N)})^2$, the eigenvalues of $M$. From the general theory of circulant matrices, the eigenvalues of $M$ are 
$$
\mu_k=\sum_{j=0}^{\frac Nd-1}c_j(e^{2d\pi i/N})^{jk} =2-e^{4dk\pi i/N}-e^{2dk(\frac Nd-2)\pi i/N}=2 - e^{4dk\pi i/N}-e^{-4dk\pi i/N}= 2 - 2\cos\left(\frac{4dk\pi}{N}\right), 
$$
for $k\in\{0,\ldots, \frac{N}{d}-1\}.$ Since $\frac{N}{d}$ is odd, multiplication by $2$ just permutes the residues modulo $\frac{N}{d}$. The closest angles to $\pi$ that could maximize $\mu_k$ are therefore $\pi \pm \pi\frac{d}{N}$, and these give
\begin{equation*}
\max_k\mu_k
=2+2\cos\left(\frac{\pi d}{N}\right)
=4\cos^2\left(\frac{\pi d}{2N}\right).
\end{equation*}
Moreover, since $d = \mathrm{gcd}(n_1, n_2, N)$, the values of $n \star m$ modulo $N$ range over the multiples of $d$. Thus,
\begin{equation*}
\max_i \, (\alpha_{n, m_i})^2 = \frac{N^2}{4\pi^2}\cos^2\left( \frac{\pi d}{2N} \right).
\end{equation*}
Combining these gives
\begin{equation*}
\|\ad_{C_n}\|_{op} = \frac{N}{\pi}\cos^2\left(\frac{\pi d}{2N}\right), \qquad \|\ad_{\Lambda C_n}\|_{op} = \frac{\| n \|^2}{\pi}\cos^2\left(\frac{\pi d}{2N}\right),
\end{equation*}
since $\Lambda C_n = (\|n\|^2/N)C_n$. 
A similar argument applies to geodesics with initial velocity $u_0=S_n$. Indeed,~\eqref{eq_lie_bracket2} implies
\begin{align*}
    -\ad_{S_n}\ad_{S_n}S_m&=(\alpha_{n,m}^{(N)})^2(2S_m+S_{m-2n}+S_{m+2n}),\\
    -\ad_{S_n}\ad_{S_n}C_m&=(\alpha_{n,m}^{(N)})^2(2C_m+C_{m-2n}+C_{m+2n}),
\end{align*}
 from which we see that the matrix for $-\ad_{S_n}^2$ is again block-diagonal with two identical blocks. Considering just the $C_m$ block and ordering this portion of the basis in the same way as before leads to a matrix representation of $-\ad_{S_n}^2$ identical in form as~\eqref{ad^2matrix}, but with 
$$M=\begin{bmatrix}
    2&0&1&0&\cdots&1&0\\
    0&2&0&1&&&1\\
    1&0&2&0&\ddots&&\vdots\\
    0&1&0&
    2&&\ddots&\vdots\\
    
    \vdots&&\ddots&&\ddots&&1\\
    1&&&\ddots&&2&0\\
    0&1&\dots&\dots&1&0&2
\end{bmatrix}.$$
The eigenvalues of this circulant matrix $M$ are
$$\lambda_k=\sum_{j=0}^{\frac Nd-1}c_j(e^{2d\pi i/N})^{jk}=2+2\cos(4dk\pi/N).$$
These attain their maximum $\lambda_k = 4$ at $k = 0$. Since the block coefficients $(\alpha_{n, m_i}^{(N)})^2$ are the same as in the $C_n$ case, we obtain
\begin{equation*}\label{eq_ad_sn}
\|\ad_{S_n}\|_{op}=\frac{N}{\pi}\cos\Big(\frac{\pi d}{2N}\Big), \qquad
\|\ad_{\Lambda S_n}\|_{op}=\frac{\|n\|^2}{\pi}\cos\Big(\frac{\pi d}{2N}\Big).
\end{equation*}
In particular, $\|\ad_{\Lambda C_n}\|_{op}$ and $\|\ad_{\Lambda S_n}\|_{op}$ differ for finite $N$, but both converge to $\|n\|^2/\pi$ as $N \rightarrow \infty$ with $n$ fixed, since $d$ divides $\mathrm{gcd}(n_1, n_2)$ and therefore remains bounded.
\end{proof}

These operator norms along with ~\eqref{Lambdaeigenvalues} yield the desired bounds for the conjugate time $\tau$ provided by~\eqref{conjugatetime}, which we summarize in the following proposition.

\begin{proposition}
Suppose $\gamma$ is a geodesic on $SU(N)$ with respect to the Zeitlin metric whose initial velocity $u_0$ is equal to $C_n$ or $S_n$ for some $n$. Let $d = \gcd(n_1, n_2, N)$, and let \begin{equation*}
\mathfrak{c}(N) :=
\begin{cases}
\cos^2\left(\frac{\pi d}{2N}\right), &
\text{if $u_0 = C_n$}, \\
\cos\left(\frac{\pi d}{2N}\right), &
\text{if $u_0 = S_n$}. \\
\end{cases}
\end{equation*}
Then the first conjugate time $\tau$ satisfies
\begin{equation*}
\frac{2\pi^2}{\mathfrak{c}(N)\|n\|^2N}\le\tau\le\frac{\pi^2(N-1)^2}{\mathfrak{c}(N)\|n\|^2N}.
\end{equation*}
\end{proposition}

\subsection{Diameter of $SU(N)$}
We will now derive an upper bound on the diameter of the Lie group $SU(N)$ under the Zeitlin metric using~\eqref{diameter}. This is accomplished by first calculating 
\begin{equation}\label{infimumkappa}\inf_{\|u_0\|_\kappa=1}\|\ad_{\Lambda u_0}\|_{op},\end{equation}
where $\|\cdot\|_\kappa$ is the norm associated with the bi-invariant metric $\kappa$. This is then used in Theorem~\ref{thm:diameter}
to bound the quantity 
\[\inf_{\|u_0\|=1}\|\ad_{\Lambda u_0}\|_{op}\]
from below, the infimum now taken with respect to the norm $\|\cdot\|$ associated with the Zeitlin metric. To compute the infimum~\eqref{infimumkappa}, we first consider the eigenvalues of the commutator operator
\begin{equation*} 
L_{u_0} : \mathfrak{su}(N)\to \mathfrak{su}(N),\qquad L_{u_0}(X) := [X, u_0].
\end{equation*}
Note that $\ad_{u_0} = \tfrac{N}{2\pi}L_{u_0}$, so working with $L_{u_0}$ just simplifies the constants. Any $u_0\in \mathfrak{su}(N)$ is skew-Hermitian and traceless. Skew-Hermitian implies that the eigenvalues are purely imaginary and that $u_0$ is normal. Normality allows the spectral theorem to be applied. Thus, there is an eigenbasis for which 
\[u_0=\operatorname{diag}(i\lambda_1,\dots, i\lambda_N)\quad\text{with}\quad \sum_{i=1}^N\lambda_i=0,\]
    with the $\lambda_i$ not necessarily distinct. We may assume without loss of generality that the given $u_0$ is expressed in this eigenbasis, since $A=UBU^{-1}$ implies $\ad_A=C_U\circ\ad_B\circ C_U^{-1}$ where $C_U(X)=UXU^{-1}$. For such a $u_0$, consider the eigenvalue problem 
    \begin{equation}\label{ad_eigenvalues}
    L_{u_0}(X)=[X,u_0]=\omega X.
    \end{equation}
    However, $\mathfrak{su}(N)\subset M_N(\mathbb{C})$ is only a real vector space (it is not closed under multiplication by $i$). Thus, to make sense of this eigenvalue problem for $\omega\in\mathbb{C},$ we must consider~\eqref{ad_eigenvalues} on the complexification
    \[\mathfrak{su}(N)_{\mathbb{C}}:=\mathfrak{su}(N)\otimes\mathbb{C}\cong \mathfrak{sl}(N,\mathbb{C})=\{A\in M_N(\mathbb{C})\,|\,\operatorname{tr} A=0\}.\]
    Any traceless matrix $A$ can be written uniquely as $A=X+iY$, where
    \begin{equation}\label{decomp}
    X=\frac{A-A^\dagger}{2},\qquad Y=\frac{A+A^\dagger}{2i}
    \end{equation}
    are both elements of $\mathfrak{su}(N)$. Thus, $\mathfrak{su}(N)_{\mathbb{C}}\cong \mathfrak{su}(N)\oplus i\,\mathfrak{su}(N)$. This is the decomposition underlying the real Zeitlin basis $\{C_n,S_m\}$. The Zeitlin matrices $J_n$ are really elements of $\mathfrak{sl}(N,\mathbb{C}),$ viewed as the complexification of $\mathfrak{su}(N).$ 
    
   One can then check that on $\mathfrak{su}(N)_{\mathbb{C}}$, the standard basis matrices $E_{jk}$ (with $j\ne k$) for $M_N(\mathbb{C})$ are eigenvectors of $L_{u_0}$ with corresponding eigenvalues $i(\lambda_k-\lambda_j)$. These account for $N^2-N$ eigenvalues, including multiplicity. The remaining eigenvectors are the basis matrices for the $N-1$ dimensional Cartan subalgebra (in this case, traceless diagonal matrices), corresponding to the eigenvalue 0. This basis consists of the vectors
   \[\{E_{11}-E_{22}, E_{11}-E_{33},\dots, E_{11}-E_{NN}\}.\]

We now have the eigenvalues and eigenvectors of $\ad_{u_0}$ on $\mathfrak{su}(N)_{\mathbb{C}}$, which we will use to obtain its eigenvalues on $\mathfrak{su}(N)$. For the off-diagonal entries, use~\eqref{decomp} to write $E_{jk}=X_{jk}+iY_{jk}$ so that the eigenvalue problem~\eqref{ad_eigenvalues} becomes
\begin{equation*}
L_{u_0}(X_{jk})+iL_{u_0}(Y_{jk})=(\lambda_j-\lambda_k)Y_{jk}+i(\lambda_k-\lambda_j)X_{jk}.
\end{equation*}
Equating real and imaginary parts, one sees that the real form of the complex eigenspace of $i(\lambda_k-\lambda_j)$ is the two-dimensional real vector space $V_{jk}:=\operatorname{span}\{X_{jk},Y_{jk}\}$ on which $L_{u_0}$ acts via the rotation
\[L_{u_0}(X_{jk},Y_{jk})=\begin{pmatrix}
    0 & \lambda_j-\lambda_k \\
     \lambda_k-\lambda_j & 0 
\end{pmatrix}\begin{pmatrix}
    X_{jk}\\
    Y_{jk}
\end{pmatrix}.
\]
This yields the direct sum decomposition
\[\mathfrak{su}(N)=\ker(\ad_{u_0}) \bigoplus_{j<k}V_{jk}, \]
with the operator norm of $\ad_{u_0}=\frac{N}{2\pi}L_{u_0}$ given by 
\begin{equation*}
\|\ad_{u_0}\|_{op} = \frac{N}{2\pi} \max_{j\ne k}|\lambda_j-\lambda_k|.
\end{equation*}
Up to the factor of $\frac{N}{2\pi}$, the infimum of $\|\ad_{u_0}\|_{op}$ over unit vectors in $\mathfrak{su}(N)$ is obtained by solving the optimization problem 
\begin{equation*}
\inf_{\|u_0\|_\kappa=1}|\lambda_{max}(u_0)-\lambda_{min}(u_0)|\qquad\text{subject to}\quad\sum_{n=1}^N\lambda_n=0,\quad\sum_{n=1}^N\lambda^2_n=1.
\end{equation*}
This is done in the following theorem, thereby obtaining the diameter of $SU(N)$ under the Zeitlin metric.

\begin{theorem}\label{thm:diameter}
    The diameter of $SU(N)$ under the Zeitlin metric satisfies the bound
\begin{equation*}
\operatorname{diam}(SU(N))\leq 2\pi^2\lambda_{max}(\Lambda)=\frac{\pi^2(N-1)^2}{N}.
\end{equation*}
In particular, the diameter grows at most linearly in $N$.
\end{theorem}

\begin{proof}
    We first show that
    \begin{equation}\label{infimumbound}
        \inf_{\|u_0\|_\kappa=1}|\lambda_{max}(u_0)-\lambda_{min}(u_0)|\ge \frac{2}{\sqrt{N}}.
    \end{equation}
    Let $a=\lambda_{min}(u_0)$ and $b=\lambda_{max}(u_0)$. Then any other eigenvalue $\lambda_n$ of $u_0$ satisfies $a\le\lambda_n\le b$. Let
    \begin{equation*}
    x_n :=\frac{1}{2}(a+b) - \lambda_n,\quad\text{so that}\quad|x_n|
    \le
    \frac{b-a}{2}.
    \end{equation*}
    The zero trace constraint implies
    \[\sum_{n=1}^N x_n=-\frac{N}{2}(a+b).\]
    This, along with the unit norm constraint, gives 
    \begin{align*}
        1=\sum_{n=1}^N\lambda^2_n
         &=\sum_{n=1}^Nx^2_n-\frac{N}{2}(a+b)^2\le\sum_{n=1}^Nx^2_n\le\sum_{n=1}^N \frac{(a-b)^2}{4}=\frac{N}{4}(a-b)^2,
    \end{align*}
    which implies 
    \[\frac{2}{\sqrt{N}}\le |a-b|.\]
    This proves~\eqref{infimumbound}. Now, define the sets
    \[U=\{v\in\mathfrak{su}(N)\,:\,1\le \|\Lambda^{-1}v\|\}\quad\text{and}\quad V=\{v\in\mathfrak{su}(N)\,:\,\sqrt{\lambda_{min}(\Lambda)}\le \|v\|_\kappa\}.\]
      Then $U\subseteq V$, since 
     \[1\le\|\Lambda^{-1}v\|^2=\kappa(\Lambda^{-1}v,v)\le\lambda_{max}(\Lambda^{-1})\kappa(v,v),\]
     and $\lambda_{max}(\Lambda^{-1})=1/\lambda_{min}(\Lambda)$. Thus we have
\begin{equation*}
\inf_{\|\Lambda^{-1}v\|\ge1}\|\ad_{v}\|_{op}
\ge
\inf_{\|v\|_\kappa\ge\sqrt{\lambda_{min}(\Lambda)}}\|\ad_{v}\|_{op}
=
\sqrt{\lambda_{min}(\Lambda)}\inf_{\|v\|_\kappa\ge1}\|\ad_{v}\|_{op}
\ge
\frac{N}{2\pi}\frac{2\sqrt{\lambda_{min}(\Lambda)}}{\sqrt{N}} = \frac{1}{\pi}
\end{equation*}
using~\eqref{infimumbound}. By~\eqref{diameter}, this immediately yields the desired upper bound. 
\end{proof}

\appendix
\section{Example: $\mathfrak{su}(3)$}

We show that the reduced geodesic equation can be solved explicitly within the subalgebra
$$\mathfrak{C}(3) = \mathrm{span}\big\{ C_{(1, 0)}, C_{(0, 1)}, C_{(1, 1)}, C_{(1, -1)}\big\} \subseteq \mathfrak{su}(3).$$
In this subalgebra, the structure constants all share a factor of  $\frac{3 \sqrt{3} }{4\pi}$, which we set to $1$ for simplicity. Let
\begin{equation*}
u(t) = u_1(t)C_{(1, 0)} + u_2(t)C_{(0, 1)} + u_3(t)C_{(1, 1)} + u_4(t)C_{(1, -1)}.
\end{equation*}
Then the reduced geodesic equation $u' + \adstar_uu = 0$ becomes the system:
\begin{equation} \label{eq_u_geodesic_su3}
\begin{alignedat}{2}
u_1' &= -u_2(u_3 - u_4), &\qquad
u_2' &= u_1(u_3 - u_4), \\ 
u_3' &= \tfrac{1}{2} u_4(u_1 - u_2), &\qquad 
u_4' &= -\tfrac{1}{2} u_3(u_1 - u_2). 
\end{alignedat}
\end{equation}
This suggests the change of variables:
\begin{equation*}
v_1 = u_1+u_2, \quad v_2 = u_1-u_2, \quad v_3 = u_3+u_4, \quad v_4 = u_3-u_4
\end{equation*}
which turns \eqref{eq_u_geodesic_su3} into the system
\begin{equation*}
v_1' = v_2v_4,
\quad v_2' = -v_1v_4,
\quad v_3' = -\tfrac{1}{2}v_2v_4,
\quad v_4' = \tfrac{1}{2}v_2v_3.
\end{equation*}
From this, it follows that the quantities
\begin{equation*}
m_1^2 := v_1^2 + v_2^2, \quad
m_2^2 := v_3^2 + v_4^2, \quad
m_3 := \tfrac{1}{2} v_1 + v_3
\end{equation*}
are conserved. In particular, solutions are confined to the torus defined by the two equations $m_1^2 = const$ and $m_2^2 = const$. To describe them more concretely, let
\begin{equation*}
v_1(t) = m_1\cos(\theta), \quad 
v_2(t) = m_1\sin(\theta), \quad 
v_3(t) = m_2\cos(\phi), \quad
v_4(t) = m_2\sin(\phi).
\end{equation*}
Then the system further reduces to two unknown functions $\theta(t)$ and $\phi(t)$ satisfying
\begin{equation*}
\theta'(t) = -m_2\sin\big(\phi(t)\big), \qquad 
\phi'(t) = \frac{m_1}{2}\sin\big(\theta(t)\big).
\end{equation*}
Since $m_3$ is also conserved, we get the further constraint that the motion $(\theta(t), \phi(t)) \in T^2$ obeys
\begin{equation*}
\frac{m_1}{2}\cos\big(\theta(t)\big) + m_2\cos\big(\phi(t)\big) = const = m_3.
\end{equation*}
Solutions of the system can be described by Jacobi elliptic functions (see \cite{abramowitz1964handbook}), and are periodic, with two notable exceptions:
\begin{enumerate}
    \item If $m_3 = \pm \big(\tfrac{m_1}{2} + m_2\big)$, the solution is constant (steady geodesic). \\
    \item If $m_3 = \pm \big(\tfrac{m_1}{2} - m_2\big)$, the solution is either constant or an asymptotic trajectory.
\end{enumerate}
Both of these exceptions correspond to special cases in which the Jacobi elliptic functions degenerate to standard trigonometric or hyperbolic functions. For example, if $m_3 = 0$, $m_1 = 2$ and $m_2 = 1$,
\begin{equation*}
\begin{alignedat}{2} 
\cos(\theta(t)) &= \tanh(t - t_0), &\qquad \sin(\theta(t)) &= \mathrm{sech}(t - t_0), \\
\cos(\phi(t)) &= -\tanh(t - t_0), &\qquad 
\sin(\phi(t)) &= \mathrm{sech}(t - t_0), 
\end{alignedat}
\end{equation*}
for some arbitrary $t_0$.

\noindent If we consider the full Lie algebra $\mathfrak{su}(3)$, then the system is
\begin{align*}
    \dot a&=-\frac{3\sqrt 3}{4\pi}(bc-bd+fg+fh)&\dot b&=-\frac{3\sqrt 3}{4\pi}(-ac+ad-eg+eh)\\
    \dot c&=-\frac{3\sqrt 3}{4\pi}\left(-\frac{1}{2}ad+\frac{1}{2}bd+\frac{1}{2}eh+\frac{1}{2}fh\right)&\dot d&=-\frac{3\sqrt 3}{4\pi}\left(\frac{1}{2}ac-\frac{1}{2}bc-\frac{1}{2}eg+\frac{1}{2}fg \right)\\
    \dot e&=-\frac{3\sqrt 3}{4\pi}(bg-bh-cf-df)&\dot f&=-\frac{3\sqrt 3}{4\pi}(-ag-ah+ce+de)\\
    \dot g&=-\frac{3\sqrt 3}{4\pi}\left(\frac12ah+\frac{1}{2}bh+\frac{1}{2}de-\frac{1}{2}df\right)&\dot h&=-\frac{3\sqrt 3}{4\pi}\left(-\frac{1}{2}ag-\frac{1}{2}bg-\frac{1}{2}ce-\frac{1}{2}cf\right)
\end{align*}
where $u(t)=(a(t),b(t),c(t),d(t),e(t),f(t),g(t),h(t)).$

\bibliographystyle{amsplain}
\bibliography{bibliography}

\end{document}